\documentclass[11pt,a4paper]{amsart}

\usepackage{amssymb}
\swapnumbers
\makeatletter
\def\section{\@startsection{section}{1}%
  \z@{.7\linespacing\@plus\linespacing}{.5\linespacing}%
  {\normalfont\bfseries\centering}}
\def\@secnumfont{\bfseries}
\makeatother
\def\frak{\mathfrak}
\def\Bbb{\mathbb}
\def\Cal{\mathcal}

\newcommand{\xbb}[3]{\begin{picture}(46,12)\put(3,3){\line(1,0){40}}%
\put(3,3){\makebox(0,0){$\times$}}\put(23,3){\makebox(0,0){$\bullet$}}%
\put(43,3){\makebox(0,0){$\bullet$}}\put(3,10){\makebox(0,0){\scriptsize $#1$}}%
\put(23,10){\makebox(0,0){\scriptsize $#2$}}%
\put(43,10){\makebox(0,0){\scriptsize $#3$}}\end{picture}}

\newcommand{\xbbdb}[5]{\begin{picture}(66,12)\put(3,3){\line(1,0){40}}%
\put(53,3){\makebox(0,0){\dots}}%
\put(3,3){\makebox(0,0){$\times$}}\put(23,3){\makebox(0,0){$\bullet$}}%
\put(43,3){\makebox(0,0){$\bullet$}}\put(63,3){\makebox(0,0){$\bullet$}}%
\put(3,10){\makebox(0,0){\scriptsize $#1$}}%
\put(23,10){\makebox(0,0){\scriptsize $#2$}}%
\put(43,10){\makebox(0,0){\scriptsize $#3$}}%
\put(53,10){\makebox(0,0){\scriptsize $#4$}}%
\put(63,10){\makebox(0,0){\scriptsize $#5$}}\end{picture}}

\newcommand{\xbbdx}[5]{\begin{picture}(66,15)\put(3,3){\line(1,0){40}}%
\put(53,3){\makebox(0,0){\dots}}%
\put(3,3){\makebox(0,0){$\times$}}\put(23,3){\makebox(0,0){$\bullet$}}%
\put(43,3){\makebox(0,0){$\bullet$}}\put(63,3){\makebox(0,0){$\times$}}%
\put(3,10){\makebox(0,0){\scriptsize $#1$}}%
\put(23,10){\makebox(0,0){\scriptsize $#2$}}%
\put(43,10){\makebox(0,0){\scriptsize $#3$}}%
\put(53,10){\makebox(0,0){\scriptsize $#4$}}%
\put(63,10){\makebox(0,0){\scriptsize $#5$}}\end{picture}}

\newcommand{\xbdbx}[5]{\begin{picture}(66,15)\put(3,3){\line(1,0){20}}%
\put(33,3){\makebox(0,0){\dots}}\put(43,3){\line(1,0){20}}%
\put(3,3){\makebox(0,0){$\times$}}\put(23,3){\makebox(0,0){$\bullet$}}%
\put(43,3){\makebox(0,0){$\bullet$}}\put(63,3){\makebox(0,0){$\times$}}%
\put(3,10){\makebox(0,0){\scriptsize $#1$}}%
\put(23,10){\makebox(0,0){\scriptsize $#2$}}%
\put(33,10){\makebox(0,0){\scriptsize $#3$}}%
\put(43,10){\makebox(0,0){\scriptsize $#4$}}%
\put(63,10){\makebox(0,0){\scriptsize $#5$}}\end{picture}}

\newcommand{\xdbbx}[5]{\begin{picture}(66,15)\put(23,3){\line(1,0){20}}%
\put(13,3){\makebox(0,0){\dots}}\put(43,3){\line(1,0){20}}%
\put(3,3){\makebox(0,0){$\times$}}\put(23,3){\makebox(0,0){$\bullet$}}%
\put(43,3){\makebox(0,0){$\bullet$}}\put(63,3){\makebox(0,0){$\times$}}%
\put(3,10){\makebox(0,0){\scriptsize $#1$}}%
\put(13,10){\makebox(0,0){\scriptsize $#2$}}%
\put(23,10){\makebox(0,0){\scriptsize $#3$}}%
\put(43,10){\makebox(0,0){\scriptsize $#4$}}%
\put(63,10){\makebox(0,0){\scriptsize $#5$}}\end{picture}}

\newtheorem*{prop*}{Proposition}

\newtheorem*{thm*}{Theorem}

\newtheorem*{lem*}{Lemma}

\newtheorem*{kor*}{Corollary}

\newcommand{\ad}{\operatorname{ad}}
\newcommand{\Ad}{\operatorname{Ad}}
\renewcommand{\exp}{\operatorname{exp}}
\newcommand{\id}{\operatorname{id}}

\renewcommand{\ker}{\operatorname{ker}}
\newcommand{\im}{\operatorname{im}}

\newcommand{\Fl}{\operatorname{Fl}}
\newcommand{\pr}{\operatorname{pr}}

\newcommand{\x}{\times}
\renewcommand{\o}{\circ}

\let\ccdot\cdot
\def\cdot{\hbox to 2.5pt{\hss$\ccdot$\hss}}

\newcommand{\al}{\alpha}

\newcommand{\ka}{\kappa}

\newcommand{\om}{\omega}
\renewcommand{\phi}{\varphi}
\newcommand{\ph}{\varphi}
\newcommand{\ps}{\psi}
\newcommand{\si}{\sigma}

\newcommand{\ze}{\zeta}

\newcommand{\Ga}{\Gamma}
\newcommand{\La}{\Lambda}
\newcommand{\Om}{\Omega}
\newcommand{\Ph}{\Phi}

\newcommand{\Up}{\Upsilon}
\newcommand{\Si}{\Sigma}

\def\sideremark#1{\ifvmode\leavevmode\fi\vadjust{
\vbox to0pt{\hbox to 0pt{\hskip\hsize\hskip1em
\vbox{\hsize3cm\tiny\raggedright\pretolerance10000
\noindent #1\hfill}\hss}\vbox to8pt{\vfil}\vss}}}

\begin{document}
\title[correspondence and twistor spaces]{Correspondence spaces and
twistor spaces for parabolic geometries}
\author{Andreas \v Cap}
\address{Institut f\"ur Mathematik, Universit\"at Wien, Nordbergstra\ss e 15,
A--1090 Wien, Austria and International Erwin Schr\"odinger Institute for
Mathematical Physics, Boltzmanngasse 9, A--1090 Wien, Austria}
\email{Andreas.Cap@esi.ac.at}

\begin{abstract}
For a semisimple Lie group $G$ with parabolic subgroups $Q\subset
P\subset G$, we associate to a parabolic geometry of type $(G,P)$ on a
smooth manifold $N$ the correspondence space $\Cal CN$, which is the
total space of a fiber bundle over $N$ with fiber a generalized flag
manifold, and construct a canonical parabolic geometry of type $(G,Q)$
on $\Cal CN$.

Conversely, for a parabolic geometry of type $(G,Q)$ on a smooth
manifold $M$, we construct a distribution corresponding to
$P$, and find the exact conditions for its integrability. If these
conditions are satisfied, then we define the twistor space $N$ as a
local leaf space of the corresponding foliation. We find equivalent
conditions for the existence of a parabolic geometry of type $(G,P)$
on the twistor space $N$ such that $M$ is locally isomorphic to the
correspondence space $\Cal CN$, thus obtaining a complete local
characterization of correspondence spaces.

We show that all these constructions preserve the subclass of normal
parabolic geometries (which are determined by some underlying
geometric structure) and that in the regular normal case, all
characterizations can be expressed in terms of the harmonic curvature
of the Cartan connection, which is easier to handle. Several examples
and applications are discussed.
\end{abstract}

\maketitle

\section{Introduction}\label{1}

This paper is devoted to the study of relations between different
geometric structures via the construction of correspondence spaces and
twistor spaces. The structures we deal with are the so--called
parabolic geometries, which may be viewed as curved analogs of
homogeneous spaces of the form $G/P$, where $G$ is a semisimple Lie
group and $P\subset G$ is a parabolic subgroup. Parabolic geometries
form a rather large class of structures, including for example
projective, conformal and non-degenerate CR--structures of hypersurface
type, as well as certain higher codimension CR structures. 

The starting point of twistor theory was R.~Penrose's idea to
associate to the Grassmannian $Gr_2(\Bbb C^4)$ of planes in $\Bbb
C^4$, which is viewed as compactified complexified Minkowski space,
the twistor space $\Bbb CP^3$, and to study the conformal geometry of
$Gr_2(\Bbb C^4)$ via the complex geometry of the twistor space. The
connection between these two manifolds is the {\em correspondence
space\/} $F_{1,2}(\Bbb C^4)$, the flag manifold of lines in planes in
$\Bbb C^4$, which canonically fibers over $Gr_2(\Bbb C^4)$ and over
$\Bbb CP^3$ and defines a correspondence between the two spaces. This
correspondence gives rise to the Penrose transform, which is a basic
element of twistor theory.

Later on, twistor theory and the Penrose transform have been extended
in two directions. On one hand, the original correspondence is
homogeneous under the simple Lie group $SL(4,\Bbb C)$, which governs
the projective geometry of $\Bbb CP^3$ and (via the isomorphism with
$SO(6,\Bbb C)$) the conformal geometry of $Gr_2(\Bbb C^4)$. The
subgroups of $SL(4,\Bbb C)$ leading to these two homogeneous spaces
are parabolic subgroups. Now one may replace $SL(4,\Bbb C)$ by an
arbitrary semisimple Lie group $G$ and fix two parabolic subgroups
$P_1,P_2\subset G$ such that $P_1\cap P_2\subset G$ is parabolic. Then
the natural fibrations from $G/(P_1\cap P_2)$ onto $G/P_1$ and $G/P_2$
define a correspondence. This is the subject of the book
\cite{Baston-Eastwood}, in which these correspondences and the
resulting Penrose transforms are studied, and applications to
differential geometry and representation theory are described.

On the other hand, it has been know for quite some time that certain
geometric structures, like conformal, projective, or almost
quaternionic structures can be viewed as ``curved analogs'' of
homogeneous spaces of the form $G/P$ as above. In the case of
four--dimensional conformal structures, twistor theory has been first
extended to curved situations in \cite{Atiyah-Hitchin-Singer} and
\cite{Penrose}. Later on, the techniques were generalized to further
geometries and they have found many applications, see for example
\cite{Merkulov} and the collection \cite{Bailey-Baston}. There are two
important remarks to be made at this point. One is that twistor theory
in many cases requires a restriction on the geometric structure, like
self duality in four--dimensional conformal geometry, or
torsion freeness in quaternionic geometry. The second important point
is that compared to the flat versions of the correspondence discussed
above, the situation loses its symmetry. While the passage ``up''
from the original manifold to the correspondence space is very similar
to the flat case, the passage ``down'' is given by passing to the leaf
space of a certain foliation, so usually this is only possible
locally. More drastically, usually there is no local geometric
structure on the twistor space, but it is only a smooth or complex
manifold.

Using E.~Cartan's concept of a generalized space, one may associate to
any homogeneous space a geometric structure defined via Cartan
connections on suitable principal bundles. In the case of homogeneous
spaces of the form $G/P$ as above the corresponding geometric
structures are called {\em parabolic geometries\/}. A general
construction of Cartan connections of this type (under small technical
restrictions) was given in the pioneering work of N.~Tanaka, see
\cite{Tanaka79}. In \cite{Morimoto}, these results were embedded into
the more general theory of Cartan connections associated to geometric
structures on {\em filtered manifolds\/}, i.e.~manifolds endowed with
a filtration of the tangent bundle by subbundles.

The main emphasis in these considerations was the solution of the
equivalence problem, while the geometry of the structures in question
was only a secondary issue. Although there were some applications of
Tanaka theory to geometrical problems which are quite close to twistor
theory, see e.g.~\cite{Takeuchi}, these works unfortunately never
became well known to people working in twistor theory, which may also
be due to their rather complicated and technical nature. 

During the last years there was a renewed interest in this class of
structures and it turned out that apart from containing interesting
examples it can be studied in a remarkably uniform way. Besides new
constructions of the Cartan connections and descriptions of the
underlying structures (see \cite{Cap-Schichl}), there is a general
theory of classes of preferred connections for parabolic geometries
(see \cite{Weyl}) which generalizes Weyl connections in conformal
geometry. This leads to a systematic way of expressing the curvature
of the Cartan connection (which is the essential invariant of any
parabolic geometry) in terms of underlying data and to powerful tools
like normal coordinates for arbitrary parabolic geometries. Finally,
there are strong general results on invariant differential operators,
i.e.~differential operators intrinsic to a parabolic geometry (see
\cite{BGG} and \cite{David-Tammo}).

The impact of this in the direction of twistor theory is that in many
cases the geometric structures showing up ``on top'' of the
correspondence (i.e. in the place of the correspondence space) are
highly interesting and much more subtle than the geometric structures
showing up ``downstairs''. In this paper, we will show that the
construction of a correspondence space works in the curved setting
without restrictions, so one may always pass ``up''. Then we will give
a complete local characterization of the geometries obtained in that
way. This is done in two steps: first we find conditions for the
existence of a twistor space, which is a candidate for a space
carrying the ``downstairs'' structure; secondly, we derive the
conditions which ensure the existence of this geometric
structure. Combining the results for going up and going down, we
obtain the precise conditions on the existence of twistor
correspondences in the classical sense.

Let us describe the contents of the paper in a little more detail. 
Technically speaking, we have to deal only with one side of the
correspondence. Thus, we consider a semisimple Lie group $G$ and two
parabolic subgroups $Q\subset P\subset G$. Starting with any parabolic
geometry of type $(G,P)$ on a manifold $N$, a simple construction leads
to the {\em correspondence space\/} $\Cal CN$, which is the total
space of a natural fiber bundle over $N$ with fiber a generalized flag
manifold, and canonically carries a parabolic geometry of type
$(G,Q)$. In fact, the Cartan connections (and thus also their
curvatures) are the same for both structures, which implies some
simple restrictions on the curvature of correspondence spaces. It is
less obvious but still rather simple, that the normalization
conditions for Cartan connections of the two different types are
compatible, so for a normal parabolic geometry the correspondence
space is normal, too.

The main result of section \ref{2} is that the simple curvature
restrictions from above actually characterize correspondence spaces
locally. If $M$ is a smooth manifold equipped with a parabolic
geometry of type $(G,Q)$, then the subalgebra $\frak p\subset\frak g$
gives rise to a distribution on $M$. Under a weak torsion condition,
this distribution is integrable, so we can consider a local leaf--space
for the corresponding foliation, which is then called the {\em twistor
space\/} $N$ of $M$. The main result is then that if $M$ satisfies the
curvature restriction for correspondence spaces from above, then there
is a parabolic geometry of type $(G,P)$ on the twistor space $N$
(which is uniquely determined provided that $P/Q$ is connected) such
that $M$ is locally isomorphic to $\Cal CN$.

To get to the classical form of the twistor correspondence, one starts
with two parabolics $P_1,P_2\subset G$ which contain the same Borel
subgroup, takes a parabolic geometry of type $(G,P_1)$, constructs the
correspondence spaces associated to $P_1\cap P_2\subset P_1$ and then
the twistor space corresponding to $P_2\supset P_1\cap P_2$. It should
be emphasized that for the existence of the Penrose transform, a
parabolic geometry on the twistor space is not needed, since one side
of the transform (the pull back part) does not need any geometric
structure. On the other side one then has a fiber bundle with fiber a
generalized flag manifold, and thus results from representation theory
apply to the push down part, which is much more subtle.

The curvature of the canonical Cartan connection is a rather
complicated object in general. A well known nice feature of regular
normal parabolic geometries is that one can pass from this curvature
to the harmonic curvature, which is a much simpler object. In section
\ref{3}, we show how rather deep general results on parabolic
geometries can be used to prove that for regular normal parabolic
geometries, the curvature restrictions from above are equivalent to
the analogous restrictions for the harmonic curvature, which are much
easier to verify. Hence we arrive at effective conditions for the
existence of a twistor space, as well as for the existence of a
parabolic geometry on this twistor space. We show that these results also
work in the holomorphic category.

In section \ref{4}, we discuss several examples and outline some
applications. We first discuss in detail the example of Lagrangian
contact structures, which leads to an interesting geometric
interpretation of the projective curvature of a linear connection, as
well as to results on contact structures and partial connections on
projectivized tangent bundles. Next, we briefly outline the case of
elliptic partially integrable almost CR structures of CR dimension and
codimension two, in which even the construction of correspondence
spaces leads to unexpected results. Finally, we briefly discuss the
example of almost Grassmannian structures which leads to an
interpretation of path geometries as parabolic geometries and contains
the classical twistor theory for split signature conformal
four--manifolds.

\subsection*{Acknowledgments:} This work was supported by
project P~15747--N05 of the Fonds zur F\"orderung der
wissenschaftlichen Forschung (FWF). Discussions with M.~Eastwood,
R.~Gover, D.~Grossman, P.~Michor, G.~Schmalz, J.~Slovak, V.~Sou\v cek,
and J.~Teichmann have been extremely important. I would also like to
thank the referee for his profound and helpful comments.

\section{Correspondence spaces and twistor spaces}\label{2} 

We start by briefly reviewing some general facts about 
parabolic geometries, see \cite{Cap-Schichl}, \cite{BGG}, and
\cite{Weyl} for more information.

\subsection{$|k|$--graded Lie algebras and parabolic geometries}\label{2.1}
Let $\frak g$ be a (real or complex) semisimple Lie algebra endowed
with a $|k|$--grading, i.e.~a grading of the form $\frak g=\frak
g_{-k}\oplus\dots\oplus\frak g_0\oplus\dots\oplus\frak g_k$ such that
$\frak g_1$ generates the subalgebra $\frak p_+:=\frak
g_1\oplus\dots\oplus\frak g_k$, and such that none of the simple
ideals of $\frak g$ is contained in $\frak g_0$. Define $\frak p$ to
be the subalgebra $\frak g_0\oplus\frak p_+$. Since we will deal with
different parabolics at the same time, we will write $\frak p_0$ for
$\frak g_0$ and we will write $\frak p_-$ for the subalgebra $\frak
g_{-k}\oplus\dots\oplus\frak g_{-1}$, which is usually called $\frak
g_-$. Be aware of the fact that $\frak p_-$ is not contained in $\frak
p$.  While the grading of $\frak g$ is not $\frak
p$--invariant, it gives rise to an invariant decreasing filtration
$\frak g=\frak g^{-k}\supset\frak g^{-k+1}\supset\dots\supset\frak g^k$
defined by $\frak g^i:=\frak g_i\oplus\dots\oplus\frak g_k$ for all
$i=-k,\dots,k$. It turns out (see \cite{Yamaguchi}), that this
filtration is completely determined by the subalgebra $\frak p$ which
is parabolic, and conversely any parabolic subalgebra gives rise to a
$|k|$--grading.

Next, let $G$ be a Lie group with Lie algebra $\frak g$, and define
subgroups $P_0\subset P\subset G$ (usually $P_0$ is denoted by $G_0$)
as the subgroups of those elements, whose adjoint actions on $\frak g$
preserve the grading respectively the filtration of $\frak g$. One
shows that Lie algebras of $P$ and $P_0$ are $\frak p$ and $\frak
p_0$, the exponential map restricts to a diffeomorphism from
$\frak p_+$ onto a normal subgroup $P_+\subset P$, and that $P$ is the
semidirect product of $P_0$ and $P_+$.

A {\em parabolic geometry\/} of type $(G,P)$ on a smooth manifold $M$
(having the same dimension as the homogeneous space $G/P$) is given by
a principal $P$--bundle $p:\Cal G\to M$ and a {\em Cartan
connection\/} $\om\in \Om^1(\Cal G,\frak g)$, i.e.~a smooth one--form
with values in $\frak g$ such that
\begin{enumerate}
\item[(1)] $\om(\ze_A)=A$ for all fundamental fields $\ze_A$, $A\in
\frak p$
\item[(2)] $(r^g)^*\om = \Ad(g^{-1})\o \om$ for all $g\in P$
\item[(3)] $\om|_{T_u\Cal G}: T_u\Cal G\to \frak g$ is a linear isomorphism for
all $u\in \Cal G$.
\end{enumerate} 
The homogeneous model for this parabolic geometry is the principal
$P$--bundle $p:G\to G/P$ together with the left Maurer--Cartan form as
a Cartan connection. A {\em morphism\/} between parabolic geometries
$(p:\Cal G\to M,\om)$ and $(p':\Cal G'\to M',\om')$ of the same type
is a principal bundle map $F:\Cal G\to\Cal G'$ such that
$F^*\om'=\om$. This compatibility of $F$ with the Cartan connections
implies that it is a local diffeomorphism.

The {\em curvature--function\/} $\ka:\Cal G\to L(\La^2\frak g,\frak
g)$ of a parabolic geometry $(p:\Cal G\to M,\om)$ is defined by
$\ka(u)(X,Y):=d\om(\om_u^{-1}(X),\om_u^{-1}(Y))+[X,Y]$ for $u\in\Cal
G$, so this exactly measures to what extent the Maurer--Cartan
equation fails to hold. The defining properties of $\om$ imply that
$\ka$ is $P$--equivariant and it vanishes if one of its entries lies
in $\frak p\subset\frak g$. Hence we will view $\ka$ as an equivariant
smooth function $\Cal G\to L(\La^2\frak g/\frak p,\frak g)$. The
values of the curvature function can be used to define various
subcategories of parabolic geometries:

First, if $\ka$ is identically zero, then the corresponding parabolic
geometry is called {\em (locally) flat\/}. General results on Cartan
connections imply that flat parabolic geometries are locally
isomorphic (as parabolic geometries) to the homogeneous model $G/P$,
see Proposition 4.12 in \cite{Cap-Schichl} for a proof in the realm of
parabolic geometries. 

Second, the parabolic geometry is called {\em torsion--free\/} if the
curvature function has the property that $\ka(u)(X,Y)\in\frak p$ for
all $u\in\Cal G$ and $X,Y\in\frak g$.

Finally, the parabolic geometry is called {\em regular\/} if for all
$u\in\Cal G$, all $i,j<0$ and all $X\in\frak g^i$ and $Y\in\frak g^j$
one has $\ka(u)(X,Y)\in\frak g^{i+j+1}$. This means that with respect
to the grading on $\frak g$, all nonzero homogeneous components of
$\ka$ are of strictly positive degree. Notice that torsion free
parabolic geometries are automatically regular, so regularity should
be viewed as a condition avoiding particularly bad types of torsion.

\subsection{Normal parabolic geometries}\label{2.2}
Parabolic geometries are mainly studied because they provide a
conceptual way to describe certain underlying geometric structures, for
example conformal, almost quaternionic, or CR structures of
hypersurface type. This underlying structure easily leads to a
principal $P_0$--bundle $\Cal G_0\to M$ for an appropriate choice of
$G$ and $P$. Using quite sophisticated procedures, one extends this
bundle to a principal $P$--bundle $\Cal G\to M$ and constructs a
canonical Cartan connection on $\Cal G$. To make this Cartan
connection unique (up to isomorphism), one has to impose an additional
normalization condition on the curvature. One then arrives at an
equivalence of categories between the underlying geometric structures
and regular normal parabolic geometries. Different versions of such
prolongation procedures can be found in \cite{Tanaka79},
\cite{Morimoto}, and \cite{Cap-Schichl}.

For general Cartan connections the problem of finding an appropriate
normalization condition is very subtle, but for parabolic geometries
Lie theory offers a uniform approach. The Killing form of $\frak g$
defines a duality between $\frak g/\frak p$ and $\frak p_+$ which is
compatible with the natural $P$--actions on both spaces (which both
come from the restriction of the adjoint action of $G$). Hence for
each $k$ we get an isomorphism $L(\La^k\frak g/\frak p,\frak
g)\cong\La^k\frak p_+\otimes \frak g$ of $P$--modules. The latter
spaces are the groups in the standard complex computing the Lie
algebra homology of $\frak p_+$ with coefficients in $\frak g$. The
differentials in this standard complex define linear maps 
$$
\partial^*:L(\La^k\frak g/\frak
p,\frak g)\to L(\La^{k-1}\frak g/\frak p,\frak g),
$$ 
which are traditionally referred to as the {\em
codifferential\/}. From the explicit formula for theses differentials
one immediately reads off that they are $P$--homomorphisms. For the
curvature function $\ka:\Cal G\to L(\La^2\frak g/\frak p,\frak g)$ of
a parabolic geometry, we can form $\partial^*\o\ka:\Cal G\to L(\frak
g/\frak p,\frak g)$ and the geometry is called {\em normal\/} if this
composition vanishes identically.

By construction, $\partial^*\o\partial^*=0$ and the resulting complex
computes the Lie algebra homology of $\frak p_+$ with coefficients in
$\frak g$. Let us denote the $k$th homology group (which is a
$P$--module by construction) by $\Bbb H^k_{\frak g}$. One easily shows
that $P_+$ acts trivially on $\Bbb H^k_{\frak g}$, so the $P$--action
is determined by the action of $P_0$. Kostant's version of the
Bott--Borel--Weyl theorem in \cite{Kostant} can be used to
algorithmically compute the cohomology groups $H^*(\frak p_+,\frak
g)$, which are well known to be dual to the homology groups, as
representations of $P_0$.

The curvature function $\ka:\Cal G\to L(\La^2\frak g/\frak p,\frak g)$
of any parabolic geometry is by construction $P$--equivariant. For a
normal parabolic geometry, we get $\partial^*\o\ka=0$ and thus we can
consider the induced function $\ka_H:\Cal G\to\Bbb H^2_{\frak g}$,
which is again $P$--equivariant. This is called the {\em harmonic
  curvature\/} of the parabolic geometry. For a regular normal
parabolic geometry the Bianchi identity implies that $\ka$ vanishes
identically provided that $\ka_H$ vanishes identically. This reduction
to the harmonic curvature is a major simplification. Equivariancy of
$\ka_H$ implies that it determines a section of the bundle $\Cal
G\x_P\Bbb H^2_{\frak g}$. But since $P_+$ acts trivially on $\Bbb
H^2_{\frak g}$ we may identify this associated bundle with $\Cal
G_0\x_{P_0}\Bbb H^2_{\frak g}$, where $\Cal G_0:=\Cal G/P_+$. This is
exactly the bundle encoding the underlying geometric structure, so the
harmonic curvature can be directly interpreted in terms of this
underlying structure, while to understand the Cartan curvature $\ka$
one needs the full Cartan bundle $\Cal G$. 

A conceptual approach to the computation and geometric interpretation
of the harmonic curvature is offered by the so--called
Weyl--structures, see \cite{Weyl}. These are global smooth
$G_0$--equivariant sections $\si:\Cal G_0\to\Cal G$ of the natural
projection $\Cal G\to\Cal G/P_+=\Cal G_0$. In \cite{Weyl} it is shown
that such sections always exist, and how the pull-back $\si^*\ka$ can
be computed in terms of tensors and connections naturally associated
to $\si$. Now $\si^*\ka$ corresponds to a $P_0$--equivariant function
$\Cal G_0\to L(\La^2\frak g/\frak p,\frak g)$ and on the latter space
one has a $P_0$--equivariant algebraic Hodge decomposition, see
\ref{3.1} below. The harmonic part of this function can be interpreted
as a $P_0$--equivariant function $\Cal G_0\to\Bbb H^2_{\frak g}$,
which exactly represents the harmonic curvature $\ka_H$. It should
also be noted that usually there are more direct interpretations for
the components of $\ka_H$ of lowest homogeneity, in particular if they
are of torsion type, i.e.~have values in $\frak p_-$. We will discuss
this in some examples in section \ref{4}.

\subsection{}\label{2.3}
Suppose that $\frak g$ is complex, and we have fixed a Cartan
subalgebra $\frak h\subset\frak g$ and a choice of positive
roots. Then any parabolic subalgebra of $\frak g$ is conjugate to a
standard parabolic subalgebra, i.e.~one that contains $\frak h$ and
all positive root spaces. Standard parabolic subalgebras in $\frak g$
are in bijective correspondence with subsets $\Si\subset\Delta_0$ of
the set of simple roots. The corresponding $|k|$--grading of $\frak g$
is then given by the $\Si$--height, i.e.~the root space corresponding
to a root $\al$ lies in $\frak g_j$, where $j$ is the sum of the
coefficients of all elements of $\Si$ in the (unique) expansion of
$\al$ as a linear combination of simple roots. In particular, $\Si$
consists of those simple roots $\al$, for which the $(-\al)$--root
space is {\em not\/} contained in the parabolic.

Let us now consider the case of two nested (standard) parabolic
subalgebras $\frak q\subset\frak p\subset\frak g$. By construction,
the subset $\tilde\Si$ associated to $\frak q$ has to contain the
subset $\Si$ associated to $\frak p$. In the language of Dynkin
diagrams with crossed nodes (see chapter 2 of \cite{Baston-Eastwood}
and \cite{Cap-Schichl}) this just means that in the complex case one
passes from $\frak p$ to $\frak q$ by changing any number of uncrossed
nodes to crossed nodes. In the real case one in addition has to take
care that the new parabolic subalgebra of the complexification comes
from the given real form, which can be read off the Satake diagram. 

In any case, we get the two decompositions $\frak g=\frak
p_-\oplus\frak p_0\oplus\frak p_+$ and $\frak g=\frak q_-\oplus\frak
q_0\oplus\frak q_+$. Since $\Si\subset\tilde\Si$ we conclude that
$\frak p_\pm\subset\frak q_\pm$, $\frak q_0\subset\frak p_0$, and that
$\frak q_-=\frak p_-\oplus(\frak p\cap\frak q_-)$ and $\frak p=(\frak
p\cap\frak q_-)\oplus\frak q$. Finally note that $\frak p\cap\frak
q_-=\frak p_0\cap\frak q_-$.

The first crucial observation is that the codifferentials
$\partial^*_{\frak p}$ and $\partial^*_{\frak q}$ corresponding to the
two parabolics are compatible. The inclusion $\frak
p_+\hookrightarrow\frak q_+$ induces an inclusion  
$$
j:L(\La^k\frak g/\frak p,\frak g)\cong \La^k\frak
p_+\otimes\frak g\to\La^k\frak q_+\otimes\frak g\cong L(\La^k\frak
g/\frak q,\frak g).
$$ 
Since the standard differentials for Lie algebra homology of $\frak
p_+$ and $\frak q_+$ are both restrictions of the standard
differential for Lie algebra homology of $\frak g$ we conclude 
\begin{prop*}
For the natural inclusion $j:L(\La^k\frak g/\frak p,\frak g)\to
L(\La^k\frak g/\frak q,\frak g)$ we have
$\partial^*_{\frak q}\o j=j\o\partial^*_{\frak p}$.
\end{prop*}

\subsection{}\label{2.4}
Let $G$ be a Lie group with Lie algebra $\frak g$ and let $\frak
q\subset\frak p\subset \frak g$ be parabolic subalgebras as
before. The explicit descriptions of $\frak p$ and $\frak q$, the
corresponding gradings of $\frak g$ in terms of root spaces, and the
definition of the subgroups $P,Q\subset G$ imply that $Q$ is a closed
subgroup of $P$, and $P_+\subset Q_+$. It can be shown that $P/Q$ is
the quotient of the semisimple part of $P_0$ by some parabolic
subgroup and thus a generalized flag manifold, see section 2.4 of
\cite{Baston-Eastwood}.

Now suppose that $(p:\Cal G\to N,\om)$ is a parabolic geometry of type
$(G,P)$. Since $\Cal G$ is a principal $P$--bundle, we may restrict
the principal action to the closed subgroup $Q\subset P$, which still
acts freely on $\Cal G$. Thus, the orbit space $\Cal CN:=\Cal G/Q$ is a
smooth manifold. 

\subsection*{Definition}
The {\em correspondence space\/} $\Cal CN$ of the parabolic geometry
 $(p:\Cal G\to N,\om)$ of type $(G,P)$ is the orbit space $\Cal G/Q$.

\begin{prop*}
The correspondence space $\Cal CN$ is the total space of a natural
fiber bundle over $N$ with fiber the generalized flag manifold
$P/Q$. It carries a natural parabolic geometry of type $(G,Q)$. The
curvature function $\ka^{\Cal CN}$ of this geometry is given by
$j\o\ka^N$, where $\ka^N$ denotes the curvature function of $N$ and
$j:L(\La^2\frak g/\frak p,\frak g)\to L(\La^2\frak g/\frak q,\frak g)$
is the natural inclusion. In particular:\newline
(1) $\ka^{\Cal CN}(X,Y)=0$ if $X\in\frak p/\frak q\subset\frak g/\frak
q$.\newline
(2) If one starts from a normal parabolic geometry on $N$, then the
induced parabolic geometry on $\Cal CN$ is normal, too. 
\end{prop*}
\begin{proof}
The orbit space $\Cal CN=\Cal G/Q$ is naturally isomorphic to the
associated bundle $\Cal G\x_P P/Q$, so it is the total space of a
natural fiber bundle with fiber $P/Q$. The natural projection
$\pi:\Cal G\to\Cal G/Q$ which maps any point to its orbit is a
principal $Q$--bundle. The Cartan connection $\om\in\Om^1(\Cal G,\frak
g)$ defining the parabolic geometry on $N$ is equivariant for the
action of the group $P\supset Q$ and reproduces the generators in
$\frak p\supset\frak q$ of fundamental vector fields. Hence it is also
a Cartan connection on the principal $Q$--bundle $\Cal G\to\Cal CN$
and we have obtained a natural parabolic geometry of type $(G,Q)$ on
$\Cal CN$.

Since the parabolic geometries on $N$ and $\Cal CN$ are given by the
same Cartan connection, $\ka^{\Cal CN}$ and $\ka^N$ coincide as
functions with values in $L(\La^2\frak g,\frak g)$ so we clearly get
$\ka^{\Cal CN}=j\o\ka^N$ as a function with values in $L(\La^2\frak
g/\frak q,\frak g)$. From this, (1) is obvious, while (2) follows from
Proposition \ref{2.3}. 
\end{proof}

\subsection*{Remark} 
For parabolic geometries $(p:\Cal G\to N,\om)$ and $(p':\Cal G'\to
N',\om')$, a morphism by definition is a homomorphism $\Ph:\Cal
G\to\Cal G'$ of principal $P$--bundles such that $\Ph^*\om'=\om$. Of
course, this implies that $\Ph$ is equivariant for the action of the
subgroup $Q\subset P$, so it also defines a morphism of parabolic
geometries of type $(G,Q)$ from $(\Cal G\to\Cal CN,\om)$ to $(\Cal
G'\to\Cal C N',\om)$. Thus forming the correspondence space is a
functorial construction.

Conversely, a morphism between the two correspondence spaces is a
$Q$--equivariant map $\Ph:\Cal G\to\Cal G'$ such that
$\Ph^*\om'=\om$. To understand the behavior of $\Ph$ with respect to
the action of the group $P\supset Q$, observe that $\om$ and $\om'$
both reproduce the generators of fundamental vector fields
corresponding to the $P$--action. From this one easily concludes that
$\Ph$ is equivariant for the action of elements of the form $\exp(A)$
for $A\in\frak p$. If $P/Q$ is connected, then these elements together
with $Q$ generate all of $P$, so $\Ph$ is automatically
$P$--equivariant. In particular we see that if $P/Q$ is connected,
then isomorphism of the correspondence spaces $\Cal CN$ and $\Cal CN'$
implies isomorphism of $N$ and $N'$, and in particular for any $N$ the
automorphism groups of the geometries on $N$ and $\Cal CN$ coincide.

\subsection{}\label{2.5}
The first step towards a converse of the construction of the
correspondence space is to observe that the subalgebra $\frak p$ gives
rise to a distribution on the base space of any parabolic geometry
$(\pi:\Cal G\to M,\om)$ of type $(G,Q)$. Since $\om$ is a Cartan
connection on $\Cal G$, the tangent bundle of $M$ is the associated
bundle $\Cal G\x_Q(\frak g/\frak q)$. Since $Q\subset P$, the subspace
$\frak p\subset\frak g$ is $Q$--invariant, and $\frak q\subset\frak
p$. Thus, $\frak p/\frak q\subset\frak g/\frak q$ is a $Q$--invariant
subspace which gives rise to a smooth subbundle $E:=\Cal G\x_Q(\frak
p/\frak q)$ of $TM$, i.e.~a smooth distribution of constant rank on
$M$. Explicitly, the subspace $E_x\subset T_xM$ can be described as
the space of those tangent vectors such that for one (or equivalently
any) point $u\in\Cal G$ such that $\pi(u)=x$ and one (or equivalently
any) lift $\tilde\xi$ of $\xi$ to $T_u\Cal G$ we have
$\om(\tilde\xi)\in\frak p$. From the construction of the
correspondence space it follows immediately that for $M=\Cal CN$ the
distribution $E$ is exactly the vertical subbundle with respect to the
projection $\pi:\Cal CN\to N$, so this distribution is globally
integrable with leaf--space $N$. Integrability of the
distribution $E$ is a rather weak condition in general:

\begin{prop*}
Let $(\pi:\Cal G\to M,\om)$ be a parabolic geometry of type $(G,Q)$
with curvature function $\ka$ and let $E\subset TM$ be the smooth
distribution corresponding to the subalgebra $\frak p\supset\frak
q$. Then the distribution $E$ is integrable if and only if for all
$u\in\Cal G$ and $X,Y\in\frak p/\frak q$, we have $\ka(u)(X,Y)\in\frak
p\subset\frak g$.
\end{prop*}
\begin{proof}
Locally, any vector field on $M$ can be lifted to a vector field on
$\Cal G$. Sections of the subbundle $E\subset TM$ correspond exactly
to sections of the subbundle $\om^{-1}(\frak p)\subset T\Cal G$. Since
the Lie bracket of two lifts is a lift of the Lie bracket of the
original fields, we see that $E$ is integrable provided that the space
of sections of $\om^{-1}(\frak p)\subset T\Cal G$ is closed under the
Lie bracket. By definition of the exterior derivative, this is
equivalent to $d\om(\xi,\eta)\in\frak p$ for all
$\xi,\eta\in\Ga(\om^{-1}(\frak p))$. By definition of the curvature
function,
$$
d\om(\xi,\eta)=\ka(\om(\xi),\om(\eta))-[\om(\xi),\om(\eta)]. 
$$
Since $\frak p\subset\frak g$ is a Lie subalgebra, we see that
$\ka(X,Y)\in\frak p$ for all $X,Y\in\frak p/\frak q$ implies
integrability of $E$. 

Conversely, if we find $X,Y\in\frak p/\frak q$ and a point $u\in\Cal G$
such that $\ka(u)(X,Y)\notin\frak p$, then choose representatives
$\tilde X,\tilde Y\in\frak p$ for $X$ and $Y$. The tangent vectors
$\om^{-1}(\tilde X)(u)$ and $\om^{-1}(\tilde Y)(u)$ have nonzero
projections to $TM$. Extending them to projectable sections of the
subbundle $\om^{-1}(\frak p)$, we obtain two sections whose Lie
bracket does not lie in $\om^{-1}(\frak p)$. Their projections on $M$
are sections of $E$, whose Lie bracket cannot be contained in $E$.
\end{proof}

\subsection{}\label{2.6}
Once the distribution $E\subset TM$ from \ref{2.5} is integrable, we
can locally define a {\em twistor space\/} for $M$ as a local leaf
space for the corresponding foliation, i.e.~a smooth manifold $N$
together with a surjective submersion $\ps$ from an open subset $U$ of
$M$ onto $N$ such that $\ker(T_x\ps)=E_x$ for all $x\in U$.  The
existence of local leaf spaces for integrable distributions immediately
follows from the local version of the Frobenius theorem (see
e.g.\cite{KMS}, theorem 3.22) by projecting on one factor of an
adapted chart. One easily shows that for two local leaf spaces
$\ps:U\to N$ and $\ps':U'\to N'$ there is a unique diffeomorphism
$\ph:\ps(U\cap U')\to\ps'(U\cap U')$ such that $\ph\o\ps=\ps'$.

At this stage the twistor space is just a smooth manifold and only
locally defined. However Penrose transforms may be used to interpret
geometric objects on the twistor space on the original manifold $M$,
which makes the construction interesting, even without getting some
local geometric structures on the twistor space. 

Our aim is to find conditions on the curvature function of the
geometry $(\pi:\Cal G\to M,\om)$ which enable us to define a parabolic
geometry of type $(G,P)$ on a sufficiently small local leaf space
$\ps:U\to N$ such that $U$ is isomorphic (as a parabolic geometry) to
an open subset in the correspondence space $\Cal CN$. This will be
done in two steps, which require different conditions on $\ka$. First
we have to construct a diffeomorphism from an appropriate open subset
of $\Cal G$ to an open subset in a principal $P$--bundle over $N$
which satisfies a certain equivariancy condition. For $X\in\frak g$
let us denote by $X^{\Cal G}\in\frak X(\Cal G)$ the vector field
characterized by $\om(X^{\Cal G})=X$. 
\begin{prop*}
Let $(\pi:\Cal G\to M,\om)$ be a parabolic geometry of type $(G,Q)$
and let $E\subset TM$ be the subbundle corresponding to $\frak p$ as
in \ref{2.5}. Suppose further that the curvature function $\ka$ has
the property that $\ka(u)(X,Y)=0$ for all $u\in\Cal G$ and
$X,Y\in\frak p/\frak q\subset\frak g/\frak q$.

Then for any sufficiently small local leaf space $\ps:U\to N$ of the
foliation defined by $E$, there is a $Q$--equivariant diffeomorphism
$\Ph$ from a $Q$--invariant open subset of the trivial principal
bundle $N\x P$ to a $Q$--invariant open subset of $\Cal G$ such that
$\ps\o\pi\o\Ph=\pr_1:N\x P\to N$ and such that for $A\in\frak p$ and
the fundamental vector field $\ze_A\in\frak X(N\x P)$ we get
$T\Ph\o\ze_A=A^\Cal G\o\Ph$.
\end{prop*}
\begin{proof}
For $A,B\in\frak p\subset\frak g$ consider $A^{\Cal G}, B^{\Cal
G}\in\frak X(\Cal G)$ and the Lie bracket $[A^{\Cal G},B^{\Cal
G}]$. Then we compute
$$
\om([A^{\Cal G},B^{\Cal G}])=-d\om(A^{\Cal G},B^{\Cal
  G})=[A,B]-\ka(A+\frak q,B+\frak q), 
$$ 
so the assumption on $\ka$ implies that $[A^{\Cal G},B^{\Cal
G}]=[A,B]^{\Cal G}$. Hence $A\mapsto A^{\Cal G}$ restricts to a Lie
algebra homomorphism $\frak p\to\frak X(\Cal G)$, i.e.~an action of
the Lie algebra $\frak p$ on $\Cal G$. By Lie's second fundamental
theorem, see pp. 47--49 and 58 of \cite{Palais}, this Lie algebra
action integrates to a local group action. In particular, for any
$u_0\in\Cal G$ there is an open neighborhood $W$ of $(u_0,e)$ in $\Cal
G\x P$ and a smooth map $F:W\to\Cal G$ such that
\begin{itemize}
\item If $(u,e)\in W$, then $F(u,e)=u$ and
  $\tfrac{d}{dt}|_{t=0}F(u,\exp(tA))=A^{\Cal G}(u)$ for all $A\in\frak
  p$. 
\item $F(F(u,g),h)=F(u,gh)$ provided that $(u,g)$, $(u,gh)$ and
  $(F(u,g),h)$ all lie in $W$. 
\end{itemize}
Now consider a local leaf space $\ps:U\to N$ which is so small that
there is a smooth section $\si:N\to\Cal G$ of $\ps\o\pi:\pi^{-1}(U)\to
N$. Possibly shrinking the leaf space further, we find an open
neighborhood $\tilde V$ of $e$ in $P$ such that for some set $W$ as
above we have $(\si(x),g)\in W$ and $(F(\si(x),g),e)\in W$ for all
$x\in N$ and all $g\in\tilde V$. Then we define $\Ph:N\x\tilde
V\to\Cal G$ by $\Ph(x,g):=F(\si(x),g)$. For $x\in N$ the tangent map
$T_{(x,e)}\Ph:T_xN\x\frak p\to T_{\si(x)}\Cal G$ is evidently given by
$(\xi,A)\mapsto T_x\si\cdot\xi+A^{\Cal G}(\si(x))$ so it is a linear
isomorphism. Possibly shrinking $U$ and $\tilde V$, we may assume that
$\Ph$ is a diffeomorphism onto an open subset $\tilde U\subset\Cal
G$. We may further assume that $\tilde V=\{\exp(X)\exp(B):X\in
V_1,B\in V_2\}$ where $V_1$ is an open neighborhood of zero in $\frak
p\cap\frak q_-$ such that $(X,h)\mapsto \exp(X)h$ is a diffeomorphism
from $V_1\x Q$ onto an open neighborhood $V$ of $Q$ in $P$ and $V_2$
is a ball around zero in $\frak q$.

For a fixed point $x\in N$, any vector tangent to $\{x\}\x\tilde V$
can be written as $\ze_A(x,g)=\tfrac{d}{dt}|_{t=0}(x,g\exp(tA))$ for
some $g\in\tilde V$ and $A\in\frak p$. For sufficiently small $t$, we
by construction have $\Ph(x,g\exp(tA))=F(\Ph(x,g),\exp(tA))$, and thus
$T_{(x,g)}\Ph$ maps this tangent vector to $A^{\Cal
G}(\Ph(x,g))$. Thus we see that $T\Ph\o \ze_A=A^{\Cal G}\o\Ph$ on
$N\x\tilde V$. Moreover, $T\Ph\o\ze_A$ always lies in $\om^{-1}(\frak
p)\subset T\Cal G$, which implies that $\tilde\Ph(\{x\}\x\tilde V)$ is
contained in one leaf of the foliation corresponding to the integrable
distribution $\om^{-1}(\frak p)\subset T\Cal G$. From \ref{2.5} we
conclude that the map $\ps\o\pi\o\Ph$ is constant on $\{x\}\x\tilde
V$, and since $\Ph(x,e)=\si(x)$ we conclude that
$\ps\o\pi\o\Ph=\pr_1:N\x\tilde V\to N$.

For $X\in V_1$ and $B\in V_2$ we have $\exp(X)\exp(tB)\in\tilde V$ for
all $t\in[0,1]$. Since $B^{\Cal G}$ is the fundamental vector field on
$\Cal G$ generated by $B\in\frak q$, the infinitesimal condition
$T\Ph\o\ze_B=B^{\Cal G}\o\Ph$ immediately implies that
$\Ph(x,\exp(X)\exp(B))=\Ph(x,\exp(X))\cdot\exp(B)$, where in the right
hand side we use the principal right action on $\Cal G$. Since $Q$
acts freely both on $N\x V$ and on $\Cal G$ we can uniquely extend
$\Ph$ to a $Q$--equivariant diffeomorphism from $N\x V$ to the
$Q$--invariant open subset $\{u\cdot g:u\in\tilde U,g\in
Q\}\subset\Cal G$. Since the family of fundamental vector fields on
$N\x P$ and the family of the vector fields $A^{\Cal G}$ on $\Cal G$
have the same equivariancy property, this extension still satisfies
$T\Ph\o\ze_A=A^{\Cal G}\o\Ph$ for all $A\in\frak p$.
\end{proof}

\subsection{}\label{2.8}
The second step in the construction is to use the diffeomorphism from
Proposition \ref{2.6} to carry over the Cartan connection to the
principal bundle $N\x P$, which needs an additional condition on the
curvature.

\begin{thm*}
Let $(\pi:\Cal G\to M,\om)$ be a parabolic geometry of type $(G,Q)$
whose curvature $\ka$ satisfies $\ka(u)(X,Y)=0$ for all $u\in\Cal G$,
$X\in\frak p$ and all $Y\in\frak g$. Let $\ps:U\to N$ be
a sufficiently small local leaf space for the integrable distribution
$E\subset TM$ corresponding to $\frak p/\frak q\subset\frak g/\frak
q$. Then:

\noindent
(1) $\om$ induces a Cartan connection $\tilde\om$ on the trivial
principal bundle $N\x P$, which is normal if $(\pi:\Cal G\to M,\om)$
is normal.

\noindent
(2) The parabolic geometry $(\pi^{-1}(U),\om|_{\pi^{-1}(U)})$ is
isomorphic to an open subset of the correspondence space $(N\x P\to
N\x(P/Q),\tilde\om)$. If $P/Q$ is connected, this condition determines
the parabolic geometry on $N$ up to isomorphism. 
\end{thm*}
\begin{proof}
(1) By Proposition \ref{2.6} there is a open neighborhood $V$ of $e$
in $P$ and a diffeomorphism $\Ph$ from $N\x V$ onto an open subset of
$\Cal G$ such that $\ps\o\pi\o\Ph=\pr_1:N\x V\to N$ and such that
$T\Ph\o\ze_A=A^{\Cal G}\o\Ph$ for all $A\in\frak p$. Hence we can form
$\Ph^*\om\in\Om^1(N\x V,\frak g)$, which restricts to a linear
isomorphism on each tangent space. Let us denote by $\rho$ the
principal right action of $P$ on $N\x P$. We can extend the values of
this form in $N\x\{e\}$ equivariantly to a $\frak g$--valued one--form
$\tilde\om$ on $N\x P$ by putting
$$
\tilde\om(x,g):=\Ad(g^{-1})\o(\Ph^*\om)(x,e)\o T\rho^{g^{-1}}. 
$$ 
By construction, this form restricts to an isomorphism on each tangent
space and satisfies $(\rho^g)^*\tilde\om=\Ad(g)^{-1}\o\tilde\om$ for
all $g\in P$. Since the vector fields $\ze_A$ and $A^{\Cal G}$ are
$\Ph$--related, their flows are $\Ph$ related. Thus $\Ph\o
\rho^{\exp(tA)}=\Fl^{A^{\Cal G}}_t\o\Ph$, whenever the left hand side
is defined. The curvature condition $\ka(A,Y)=0$ for all $Y\in\frak g$
and $A\in\frak p$ reads as $-\ad(A)\o\om=i_{A^{\Cal G}}d\om=\Cal
L_{A^{\Cal G}}\om$, where in the last step we have used that
$\om(A^{\Cal G})$ is constant. This infinitesimal equivariancy
condition easily implies local equivariancy, i.e.~that $(\Fl^{A^{\Cal
G}}_t)^*\om=\Ad(\exp(-tA))\o\om$, whenever the flow is defined. Hence
for $A\in\frak p$ such that $\exp(tA)\in V$ for all $t\in [0,1]$ we
have $(\rho^{\exp(A)})^*\Ph^*\om(x,e)=\Ad(\exp(A)^{-1})\o\om(x,e)$ for
all $x\in N$. This shows that $\tilde\om$ coincides with $\Ph^*\om$ on
an open neighborhood of $N$ in $N\x P$, so $\tilde\om$ is smooth on
this neighborhood and hence by equivariancy on all of $N\x P$.

By construction, $\Ph^*\om$ reproduces the generators of fundamental
vector fields. Hence $\tilde\om$ has the same property in points of
the form $(x,e)$ with $x\in N$, and thus by equivariancy in all points
of $N\x P$. We have therefore verified that $\tilde\om\in\Om^1(N\x
P,\frak g)$ is a Cartan connection, and thus defines a parabolic
geometry of type $(G,P)$ on $N$. Concerning normality, we have
observed above that $\tilde\om$ coincides with $\Ph^*\om$ on an open
neighborhood of $N\x\{e\}$. On this neighborhood, the curvature
function $\tilde\ka$ of $\tilde\om$ is given by $j\o
\tilde\ka=\ka\o\Ph$, where $j$ is the inclusion from \ref{2.3}.
Now the claim on normality follows from Proposition \ref{2.3},
equivariancy of $\tilde\ka$, and the fact that $\partial^*_{\frak p}$
is $P$--equivariant.

\noindent
(2) Take the map $\Ph:N\x V\to \Cal G$ from (1). By Proposition
\ref{2.6} we may assume $N\x V$ to be a $Q$--invariant subset of $N\x
P$ and $\Ph$ to be a $Q$--equivariant diffeomorphism onto a
$Q$--invariant open subset $\tilde U\subset\Cal G$. From (1) we know
that $\Ph^*\om=\tilde\om$ on an open subset of $N$, which has to be
$Q$--invariant since both $\Ph^*\om$ and $\tilde\om$ are
$Q$--equivariant. But this exactly means that $\Ph$ defines an
isomorphism of parabolic geometries of type $(G,Q)$ from an open
subset in $\Cal CN=N\x V/Q$ to the open subset $\pi(\tilde U)\subset
M$. Finally, we have seen in Remark \ref{2.4} that isomorphism
of the correspondence spaces implies isomorphism of the underlying
spaces provided that $P/Q$ is connected. 
\end{proof}

\begin{kor*}
Suppose that $P/Q$ is connected and let $(\Cal G\to M,\om)$ be a
parabolic geometry of type $(G,Q)$ satisfying the curvature
restriction of the theorem. Suppose that $\ps:U\to N$ is any local
leaf space for the foliation corresponding to $E$. Then $N$ carries a
canonical parabolic geometry of type $(G,P)$, which is normal if
$(\Cal G,\om)$ is normal and such that $U$ is isomorphic to an open
subset of $\Cal CN$ as a parabolic geometry.
\end{kor*}
\begin{proof}
By part (1) of the theorem we get appropriate parabolic geometries on
sufficiently small open subsets of $N$. By part (2) these locally
defined structures fit together to define a principal bundle and a
Cartan connection on $N$. Also the isomorphisms between open subsets
of the correspondence spaces and appropriate subsets of $U$ piece
together smoothly by part (2) and Remark \ref{2.4}.
\end{proof}

\subsection*{Remarks}
(1) For many structures, there are intermediate curvature conditions
lying between the one in Proposition \ref{2.5}, which ensures
existence of a twistor space, and the one in the theorem above, which
ensures existence of a parabolic geometry on the twistor space. In
some cases, these conditions imply the existence of geometric
structures on the twistor space which are weaker than a parabolic
geometry, see \ref{4.6} for an example.

\noindent
(2) To get the classical forms of the twistor correspondence, one has
to combine the constructions of cor\-res\-pon\-dence and twistor
spaces. One starts with two parabolic subalgebras $\frak p_1,\frak
p_2\subset\frak g$ which contain the same Borel subalgebra. Then
$\frak q=\frak p_1\cap\frak p_2$ is a parabolic subalgebra,
too. Starting from a parabolic geometry of type $(G,P_1)$ on a
manifold $M$ one forms the correspondence space $\Cal CM$, which
carries a parabolic geometry of type $(G,Q)$. The parabolic $\frak
p_2\supset\frak q$ gives rise to a distribution $E$ on $\Cal
CM$. Since the curvature function of $\Cal CM$ coincides with the
curvature function of $M$, the condition for integrability of this
distribution from Proposition \ref{2.5} can be expressed in terms of
the curvature of $M$ (and with the help of the results in the next
section even in terms of the harmonic curvature). Given the
corresponding curvature restrictions, one can form local twistor
spaces. A priori, they are only defined locally on $\Cal CM$, but
using the fact that the fibers $P_1/Q$ of $\Cal CM\to M$ are compact,
one obtains a full correspondence over sufficiently small open subsets
of $M$. The fibers of $\Cal CM\to M$ then descend to submanifolds in
such twistor space, and the global structure of these submanifolds
encodes the local geometry of $M$. On the other hand, under further
restrictions on the curvature, the local twistor spaces inherit local
geometric structures. Again, this will be discussed in the examples in
section \ref{4}.

\section{Reduction to harmonic curvature components}\label{3}

The curvature of a Cartan connection is a rather complicated object,
and one needs the full Cartan bundle to interpret it geometrically. As
we have noted in \ref{2.2} in the case of regular normal parabolic
geometries there is the harmonic curvature, which can be directly
interpreted in terms of the underlying structure determining the
parabolic geometry. As noted in \ref{2.2} the harmonic curvature is
still a complete obstruction to local flatness in the regular normal
case. Our aim in this section is to express (in the regular normal
case) the curvature condition from \ref{2.5}--\ref{2.8} equivalently
in terms of the harmonic curvature, which leads to much more effective
conditions. 

The key to this is that one may recover the full curvature from the
harmonic curvature by applying an invariant differential operator,
which was first shown by D.~Calderbank and T.~Diemer in
\cite{David-Tammo}. The authors suggested to use this as a
strengthening of the Bianchi identity. The relation between curvature
and harmonic curvature follows from the theory of
Bernstein--Gelfand--Gelfand (BGG) sequences. The main result of this
section is Theorem \ref{3.2}, which partly works in the realm of
general BGG sequences.

Throughout this section, we fix two parabolic subgroups $Q\subset P$
of a semisimple Lie group $G$ corresponding to subalgebras $\frak
q\subset\frak p\subset\frak g$ as in \ref{2.3}.

\subsection{BGG sequences}\label{3.1}
Any representation of the Lie group $Q$ gives rise to a natural vector
bundle on the category of parabolic geometries of type $(G,Q)$ by
forming associated bundles to the principal Cartan bundles. The first
step towards BGG sequences is the observation that restrictions of
representations of $G$ play a special role. If $\Bbb V$ is
representation of $G$, then we can also view it as a representation of
$Q\subset G$, and the corresponding natural vector bundle is called a
\textit{tractor bundle}. The main feature of these bundles is that for
each parabolic geometry $(p:\Cal G\to M,\om)$, the Cartan connection
$\om$ induces a linear connection on the tractor bundle $VM:=\Cal
G\x_Q\Bbb V$, see \cite{Cap-Gover}.  To deal with the curvature, we
have to consider $\Bbb V=\frak g$, and the corresponding tractor
bundle $\Cal AM:=\Cal G\x_Q\frak g$ is called the {\em adjoint tractor
bundle\/}.

There are two ways to extend the natural linear connection on a
tractor bundle $VM$ to an operator on $VM$--valued differential forms,
see section 2 of \cite{BGG}. First, there is the usual
\textit{covariant exterior derivative}, which we denote by
$d^\om:\Om^k(M,VM)\to\Om^{k+1}(M,VM)$. Second, one may start from the
covariant exterior derivative on the homogeneous model $G/Q$, and
extend this in a different way to all parabolic geometries of type
$(G,Q)$. This leads to the \textit{twisted exterior derivative}
$d_{\Bbb V}$. The relation between these two operators can be easily
described explicitly. The curvature function $\ka$ has values in
$L(\La^2\frak g/\frak q,\frak g)$. Projecting in the last factor to
$\frak g/\frak q$ we obtain a $Q$--equivariant map $\Cal G\to
L(\La^2\frak g/\frak q,\frak g/\frak q)$, which corresponds to a
section $\ka_-$ of the associated bundle $\La^2T^*M\otimes TM$. This
is exactly the torsion of the Cartan connection $\om$. This leads to
an insertion operator $i_\ka$, i.e.~for $\ph\in\Om^k(M,VM)$, the form
$i_\ka\ph\in\Om^{k+1}(M,VM)$ is the alternation of
$(\xi_0,\dots,\xi_k)\mapsto
\ph(\ka_-(\xi_0,\xi_1),\xi_2,\dots,\xi_k)$.  Then $d^\om\ph=d_{\Bbb
V}\ph+i_\ka\ph$, so in particular, the two operators coincide for
torsion free parabolic geometries. To unify the presentation, let us
write $d$ for either $d^\om$ or $d_{\Bbb V}$ in the sequel. 

The second step is to compress either of these sequences of first
order operators defined on $VM$--valued differential forms to a
sequence of higher order operators defined on smaller bundles. The
bundle $\La^kT^*M\otimes VM$ corresponds to the representation
$L(\La^k\frak g/\frak q,\Bbb V)$ of $Q$. Similarly as in \ref{2.3}, via
the duality $(\frak g/\frak q)^*\cong \frak q_+$ the differentials in
the standard complex computing Lie algebra homology of $\frak q_+$
with coefficients in $\Bbb V$ define $Q$--homomorphisms
$\partial^*:L(\La^k\frak g/\frak q,\Bbb V)\to L(\La^{k-1}\frak g/\frak
q,\Bbb V)$ such that $\partial^*\o\partial^*=0$. The quotient
$\ker(\partial^*)/\im(\partial^*)$ is a representation of $Q$, which
we denote by $\Bbb H^k_{\Bbb V}$. It turns out that $Q_+\subset Q$
acts trivially on this quotient, so this representation is completely
reducible and determined by the $Q_0$--action.

The $Q$--homomorphisms $\partial^*$ induce vector bundle maps
$\La^kT^*M\otimes VM\to \La^{k-1}T^*M\otimes VM$ and correspondingly
algebraic operators on $VM$--valued forms, which we all denote by the
same symbol. The kernel an image of $\partial^*$ are natural
subbundles in $\La^kT^*M\otimes VM$ and their quotient can be
naturally identified with $H^k_{\Bbb V}M:=\Cal G\x_Q\Bbb H^k_{\Bbb
V}$. Since $Q_+$ acts trivially on $\Bbb H^k_{\Bbb V}$ we can also
view $H^k_{\Bbb V}M$ as $\Cal G_0\x_{Q_0}\Bbb H^k_{\Bbb V}$, so this
admits a direct interpretation in terms of underlying structures.

The key step for the construction of BGG sequence is to construct an
invariant differential operator $S:\Ga(H^k_{\Bbb V}M)\to\Om^k(M,VM)$,
called the splitting operator, such that for all $s\in\Ga(H^k_{\Bbb
  V}M)$ we have 
\begin{itemize}
\item $\partial^*(S(s))=0$ and $\pi_H(S(s))=s$. 
\item $\partial^*(dS(s))=0$. 
\end{itemize}
Here $\pi_H$ denotes the natural algebraic projection from sections of
$\ker(\partial^*)\subset\La^kT^*M\otimes VM$ to sections of $H^k_{\Bbb
  V}M$. 

These operators were first constructed in \cite{BGG} in terms of
homomorphisms of semi--holonomic jet modules. In \cite{David-Tammo}
the authors gave a simpler construction and extended them to operators
defined on $\Om^k(M,VM)$. We will sketch this construction
next. The operator $\square^R:=d\o\partial^*+\partial^*\o d$ is an
invariant operator on $\Om^k(M,VM)$ for each $k$. Since
$\partial^*\o\partial^*=0$, this operator preserves the subspace of
sections of the bundle $\im(\partial^*)$. One shows that the
restriction of $\square^R$ to sections of this subbundle is invertible
and the inverse is again a differential operator, which then is
natural by construction. Let $L:\Om^k(M,VM)\to\Om^{k-1}(M,VM)$ be the
composition of this inverse with $\partial^*$. Given a section
$\si\in\Ga(H^k_{\Bbb V}M)$ one chooses a section $\ph$ of
$\ker(\partial^*)\subset\La^kT^*M\otimes VM$ such that
$\pi_H(\ph)=\si$ and shows that $S(\si):=\ph-Ld\ph$ is independent of
the choice of $\ph$ and this recovers the splitting operator $S$. From
this point of view it is also easy to see that $S$ is uniquely
determined by the two properties listed above. Suppose that
$\ph\in\Om^k(M,VM)$ is such that $\partial^*(\ph)=0$ and
$\partial^*d\ph=0$. Then one easily shows that $\ph-S(\pi_H(\ph))$
must be a section of $\im(\partial^*)$ which lies in the kernel of
$\square^R$, and thus vanishes identically.

In particular, we may apply this to the special case $\Bbb V=\frak g$,
$k=2$ and the curvature $\ka$. Equivariancy of the curvature function
implies that $\ka$ can be viewed as an element of $\Om^2(M,\Cal
AM)$. Normality of the parabolic geometry exactly means
$\partial^*\ka=0$. Moreover, $\ka$ coincides with the curvature of the
natural linear connection on $\Cal AM$, so the Bianchi identity for
linear connections implies $d^\om(\ka)=0$. In particular,
$\partial^*(d^\om(\ka))=0$, which implies that $\ka=S(\ka_H)$, where
$S$ is the splitting operator for $d=d^\om$.

One of the tricky points of the theory is that there seems to be no
way to write down the operator $L$ (or $S$) in a manifestly invariant
way. To get a formula, one first has to choose a Weyl structure $\si$,
see \ref{2.2}. This reduces the structure group of all bundles in
question to $Q_0$. Now the $Q_0$--submodule $\frak q_-\subset\frak g$
is complementary to $\frak q$, so $\frak g/\frak q\cong\frak q_-$ as a
$Q_0$--module. Consequently, on the level of $Q_0$, we can view the
spaces $L(\La^k\frak g/\frak q,\Bbb V)$ as the chain groups in the
standard complex computing the Lie algebra cohomology $H^*(\frak
q_-,\Bbb V)$. The differential $\partial$ in this complex is a
$Q_0$--homomorphism (but not a $Q$--homomorphism). Then
$\square:=\partial\o\partial^*+\partial^*\o\partial$ is a
$Q_0$--homomorphism and thus, after a choice of a Weyl structure,
gives rise to a bundle map on each $\La^kT^*M\otimes VM$. It turns out
$\partial$ and $\partial^*$ are adjoint with respect to a certain
inner product. The resulting algebraic Hodge decomposition shows that
as a $Q_0$--representation we have $\Bbb H^k_{\Bbb
V}\cong\ker(\square)\subset L(\La^k\frak q_-,\Bbb V)$, and this
representation is computable via Kostant's version of the
Bott--Borel--Weil theorem, see \cite{Kostant}. The Hodge decomposition
also implies that $\square$ is invertible on $\im(\partial^*)$, and we
can finally state the formula:
$$
L=\left(\sum_{i=0}^\infty(-1)^i(\square^{-1}\square^R-\id)^i\right)
\square^{-1}\partial^*. 
$$
The sum in this formula is actually finite since
$\square^{-1}\square^R-\id$ increases homogeneous degrees with respect
to a natural grading on $\Bbb V$. 

\subsection{}\label{3.2}
We next have to study algebraic properties of the operators $L$ and
$S$. Consider a $Q$--submodule $\Bbb E\subset\ker(\partial^*)\subset
L(\La^k\frak g/\frak q,\Bbb V)$. Then we get a $Q_0$--submodule $\Bbb
E_0:=\Bbb E\cap\ker(\square)\subset\Bbb H^k_{\Bbb V}$ by identifying
$L(\La^k\frak g/\frak q,\Bbb V)$ with $L(\La^k\frak q_-,\Bbb V)$ as
above. We want to find conditions which make sure that the splitting
operator $S$ maps sections of $\Cal G_0\x_{Q_0}\Bbb E_0$ to sections
of $\Cal G\x_Q\Bbb E$. For the operator corresponding to $d=d_{\Bbb
V}$ this is rather easy and has been discussed in \cite{BGG}, but for
$d=d^\om$, the problem is much more subtle.

\subsection*{Definition}
Let $\Bbb E\subset L(\La^k\frak g/\frak q,\Bbb V)$ and $\Bbb F\subset
L(\La^2\frak g/\frak q,\frak g)$ be $Q$--submodules. We say that $\Bbb
E$ is \textit{stable under $\Bbb F$--insertions} if for $\ph\in\Bbb E$
and $\ps\in\Bbb F$ we have $\partial^*(i_\ps\ph)\in\Bbb E$, where
$i_\ps\ph$ is the alternation of the map $(X_0,\dots,X_k)\mapsto
\ph(\ps(X_0,X_1)+\frak q,X_2,\dots,X_k)$ for $X_i\in\frak g/\frak q$.

\begin{thm*}
Let $(p:\Cal G\to M,\om)$ be a regular normal parabolic geometry, with
curvature $\ka$, and harmonic curvature $\ka_H$. Let $\Bbb V$ be a
representation of $G$, and let $\Bbb E\subset \ker(\partial^*)\subset
L(\La^k\frak g/\frak q,\Bbb V)$ and $\Bbb F\subset
\ker(\partial^*)\subset L(\La^2\frak g/\frak p,\frak g)$ be
$Q$--submodules. Put $\Bbb E_0:=\Bbb E\cap\ker(\square)\subset\Bbb
H^k_{\Bbb V}$ and $\Bbb F_0:=\Bbb E\cap\ker(\square)\subset\Bbb
H^2_{\frak g}$. Let $E_0M=\Cal G_0\x_{Q_0}\Bbb E_0\subset H^k_{\Bbb
V}M$ and $EM=\Cal G\x_Q\Bbb E\subset\La^kT^*M\otimes VM$ be the
corresponding subbundles, and similarly for $F_0M$ and $FM$.

\noindent
(1) The splitting operator $S:\Ga(H^k_{\Bbb V}M)\to\Om^k(M,VM)$ for
    $d=d_{\Bbb V}$ maps sections of $E_0M$ to sections of $EM$.

\noindent
(2) Suppose that $\square(\Bbb E)\subset\Bbb E$, $\Bbb E$ is stable
    under $\Bbb F$ insertions and $\ka$ has values in $FM$. Then the
    splitting operator for $d=d^\om$ maps sections of $E_0M$ to
    sections of $EM$.

\noindent
(3) If $\square(\Bbb F)\subset\Bbb F$, $\Bbb F$ is stable under $\Bbb
    F$ insertions and $\ka_H$ is a section of $F_0M$, then $\ka$ is a
    section of $FM$.
\end{thm*}
\begin{proof}
We first claim that if $\square(\Bbb E)\subset\Bbb E$, then
$\partial^*\o d_{\Bbb V}$ maps sections of $EM$ to sections of $EM$:
Let us compute in jet--modules as it is done in \cite{BGG}. So
consider the first jet prolongation $\Cal J^1(\La^k(\frak g/\frak
q)^*\otimes\Bbb V)$, which as a $Q_0$--module is isomorphic to
$(\La^k(\frak q_-)^*\otimes\Bbb V)\oplus(\frak q_-^*\otimes
\La^k(\frak q_-)^*\otimes\Bbb V)$. By lemma 2.1 of \cite{BGG} $d_{\Bbb
V}$ is induced by the homomorhism $\Cal J^1(\La^k(\frak
q_-)^*\otimes\Bbb V)\to \La^{k+1}(\frak q_-)^*\otimes\Bbb V$, which is
given by $(e,Z\otimes f)\mapsto \partial(e)+(n+1)Z\wedge f$ for
$e,f\in \La^k(\frak q_-)^*\otimes\Bbb V$ and $Z\in \frak q_-^*$. Thus,
$\partial^*\o d_{\Bbb V}$ corresponds to the homomorphism $(e,Z\otimes
f)\mapsto\partial^*\partial(e)+(n+1)\partial^*(Z\wedge f)$. Since
$\Bbb E\subset\ker(\partial^*)$ the first summand coincides with
$\square(e)$ and by formula (1.2) of \cite{BGG} the second summand
gives $-(n+1)Z\cdot f$, so the claim follows.

\noindent
(1) Let $\tilde{\Bbb E}$ be the $P$--submodule of $L(\La^k\frak
g/\frak q,\Bbb V)$ generated by $\Bbb E_0$, and let $\tilde EM$ be the
corresponding bundle, so $\tilde EM\subset EM$. By claim 2.2 of
\cite{BGG} we get $\square(\tilde{\Bbb E})\subset\tilde{\Bbb E}$, so
by the above claim $\partial^*\o d_\Bbb V$ maps sections of $E_0M$ to
sections of $\tilde EM$. By definition
$$
L\o d_{\Bbb V}=\left(\sum_{i=0}^\infty(-1)^i(\square^{-1}\square^R-\id)^i
\right) \square^{-1}\partial^*d_{\Bbb V}. 
$$

Since $\square$ preserves $\tilde{\Bbb E}$, the corresponding
algebraic operator preserves sections of $\tilde EM$ for any choice of Weyl
structure. Hence the operator $\square^{-1}$ preserves sections of the
subbundle $\im(\partial^*)\cap\tilde EM$. Finally, on the image of
$\partial^*$ the operator $\square^R$ by definition coincides with
$\partial^*\o d_{\Bbb V}$, so it preserves sections of $
\im(\partial^*)\cap\tilde EM$. Hence the splitting operator maps
sections of $E_0M$ to sections of $\tilde EM\subset EM$. 

\noindent
(2) Since $\Bbb E$ is stable under $\Bbb F$ insertions and $\ka$ has
values in $FM$ we see that $\partial^*\o i_{\ka}$ maps sections of
$EM$ to sections of $EM$. Since $\square(\Bbb E)\subset\Bbb E$, this
together with the claim implies that $\partial^*\o d^\om$ maps
sections of $EM$ to sections of $EM$. Now the result follows exactly
as in (1).

\noindent
(3) Let us denote by $L$ the operator corresponding to
$d^\om$. Further, put $\square^R=\partial^*\o d_{\Bbb V}+d_{\Bbb
  V}\o\partial^*$ and $\tilde{\square}^R=\partial^*\o
d^\om+d^\om\o\partial^*$. Since $d^\om=d_{\Bbb V}+i_\ka$, we conclude
that on the image of $\partial^*$, we have
\begin{equation*}
\square^{-1}\tilde{\square}^R-\id=\square^{-1}\square^R-\id+
\square^{-1}\partial^*i_\ka.
\tag{$*$}
\end{equation*}
From \ref{3.1} we know that $\ka=\ka_H-L(d^\om\ka_H)$. Since we deal
with a regular parabolic geometry, all nonzero homogeneous components
of $\ka$ have degree bounded from below by some $\ell>0$, and by the
Bianchi identity (see \cite{Cap-Schichl}, 4.9) the lowest nonzero
homogeneous component of $\ka$ is harmonic. Hence $\ka$ is congruent
to $\ka_H\in\Ga(FM)$ modulo elements of homogeneous degree
$\geq\ell+1$. Hence $\partial^*(i_{\ka}\ka_H)$ is congruent modulo
elements homogeneous of degree $\geq 2\ell+1$ to
$\partial^*(i_{\ka_H}\ka_H)$, and the latter element lies in
$\Ga(FM)$, since $\Bbb F$ is stable under $\Bbb F$--insertions. Hence
we conclude that $\partial^*d^\om\ka_H$ is congruent to a section of
$FM$ modulo elements of homogeneous degree $\geq 2\ell+1$. From the
definition of $L$ and formula ($*$) above, we then conclude that
$\ka=\ka_H-L(d^\om\ka_H)$ is congruent to a section of $FM$ modulo
elements of that homogeneity.

As above, this implies that $\partial^*(i_{\ka}\ka_H)$ is congruent to
a section of $FM$ modulo elements homogeneous of degree $\geq
3\ell+1$. Hence  $\partial^*d^\om\ka_H$ and therefore $\ka$ are
congruent to sections of $FM$ modulo elements of that homogeneity, and
iterating this argument, the result follows. 
\end{proof}

\begin{kor*}
Let $\Bbb E\subset\ker(\partial^*)\subset L(\La^2\frak g/\frak q,\frak
g)$ be a $Q$--submodule and put $\Bbb E_0:=\Bbb
E\cap\ker(\square)\subset \Bbb H^2_\frak g$. Let $(p:\Cal G\to M,\om)$
be a regular normal parabolic geometry such that the harmonic
curvature $\ka_H$ has values in $\Bbb E_0$. If either $(\Cal G,\om)$
is torsion free or $\square(\Bbb E)\subset\Bbb E$ and $\Bbb E$ is
stable under $\Bbb E$--insertions, then the curvature function $\ka$
has values in $\Bbb E$.
\end{kor*}

\subsection{}\label{3.3}
Suppose that $G$ is a semisimple Lie group, $\frak
q\subset\frak p\subset\frak g$ are two parabolic subalgebras as in
\ref{2.3} and $Q\subset P\subset G$ are the corresponding
subgroups. Consider a regular normal parabolic geometry $(p:\Cal G\to
M,\om)$ of type $(G,Q)$. Let $\ka$ be the curvature function, $\ka_H$
the harmonic curvature and $\ka_-$ the torsion, which we consider as a
section of the bundle $L(\La^2TM,TM)$. Let $E\subset TM$ be the
distribution corresponding to $\frak p/\frak q\subset\frak g/\frak
q$. Then we have:

\begin{thm*}
(1) The distribution $E$ is integrable if and
only if $\ka_-(\xi,\eta)\in E$ for all $\xi,\eta\in E$.

\noindent
(2) The parabolic geometry $(\Cal G,\om)$ satisfies the curvature
condition of Theorem \ref{2.8} if and only if $\ka_H$ has values in
the space of those maps which vanish if one of their entries is from
$\frak p\cap\frak q_-$. 
\end{thm*}
\begin{proof}
(1) is only a straightforward reformulation of Proposition \ref{2.5}. 

\noindent
(2) Consider the submodule $\Bbb E\subset L(\La^2\frak g/\frak q,\frak
    g)$ of those maps which vanish if one of their entries lies in
    $\frak p/\frak q$. To prove that $\ka$ has values in $\Bbb E$, by
    Corollary \ref{3.2} we have to show that $\square(\Bbb
    E)\subset\Bbb E$ and $\Bbb E$ is stable under $\Bbb
    E$--insertions. By definition $\Bbb E$ is exactly the image of the
    natural inclusion 
$$ 
j:L(\La^2\frak g/\frak p,\frak g)\cong\La^2\frak p_+\otimes\frak g
\to\La^2\frak q_+\otimes\frak g\cong L(\La^2\frak g/\frak q,\frak g)
$$ 
from \ref{2.3}. By definition of the insertion operator, for
$j(\ph),j(\ps)\in\Bbb E$, we get $i_{j(\ph)}j(\ps)=j(i_\ph\ps)\in\Bbb
E$. But by Proposition \ref{2.3} we get $\partial^*_\frak
q(j(i_\ph\ps))=j(\partial^*_{\frak p}(i_\ph\ps))\in\Bbb E$, so $\Bbb
E$ is stable under $\Bbb E$--insertions. 

To prove that $\square(\Bbb E)\subset\Bbb E$, we use Kostant's
algebraic formula for $\square$. Let $\tilde j:\La^2\frak
q_+\otimes\frak g\to \La^2\frak g\otimes\frak g$ be the
inclusion. Recall from \ref{2.3} that we have the decomposition $\frak
g=\frak p_-\oplus (\frak p\cap\frak q_-)\oplus\frak q_0\oplus (\frak
p_0\cap\frak q_+)\oplus\frak p_+$.  The Killing form of $\frak g$
induces dualities between the first and last, and the second and
fourth summands and its restriction to the third summand is
non--degenerate. Now we choose a basis $\{X_\al\}$ of $\frak g$
consisting of homogeneous elements, which starts with bases of the
first three summands and has the duals of the first two bases in the
last two summands. Denoting by $\{Y_\al\}$ the dual basis with respect
to the Killing form, we by construction have $[X_\al,Y_\al]\in\frak
p_0$ for all $\al$ and the first elements of $\{Y_\al\}$ (which lie in
$\frak p_+$) coincide with the last elements of $\{X_\al\}$ and vice
versa. According to Theorem 4.4 of \cite{Kostant}, one can write
$\tilde j\o\square$ as
$$ 
\tfrac{1}{2}\bigg(\sum_{\al}\id\otimes(\ad_{Y_{\al}}\o\ad_{X_{\al}})+
\sum_{\al:X_{\al}\in\frak q_-}\rho_{Y_{\al}}\o\rho_{X_{\al}}-
\sum_{\al:X_{\al}\in\frak q}\rho_{Y_{\al}}\o\rho_{X_{\al}}\bigg)\o \tilde
j.
$$ 
Here $\ad$ denotes the adjoint action of $\frak g$ on $\frak g$,
while $\rho$ denotes the natural action of $\frak g$ on $\La^2\frak
g\otimes\frak g$. The subspace $\im(\tilde j\o j)$ is invariant under
any map which acts only on the $\frak g$ part, as well as under
$\rho_A$ for each $A\in\frak p$. In particular, any summand in the
first sum preserves this subspace, while for the other two sums the
same holds for summands in which $X_\al$ lies in $(\frak p\cap\frak
q_-)\oplus\frak q_0\oplus(\frak p_0\cap\frak q_+)\subset\frak p$,since
then also $Y_\al$ lies in this subspace. It remains to show that
$\sum_{\al:X_{\al}\in\frak p_-}\rho_{Y_{\al}}\o\rho_{X_{\al}}-
\sum_{\al:X_{\al}\in\frak p_+}\rho_{Y_{\al}}\o\rho_{X_{\al}}$ also
preserves this subspace. But by construction, we can rewrite this part
as
$$
\sum_{\al:X_{\al}\in\frak p_-}(\rho_{Y_{\al}}\o\rho_{X_{\al}}-
\rho_{X_{\al}}\o\rho_{Y_{\al}})=\sum_{\al:X_{\al}\in\frak
  p_-}\rho_{[Y_\al,X_\al]},
$$
and since $[Y_\al,X_\al]\in\frak p$, the result follows. 
\end{proof}

\subsection{The relation between real and complex parabolic
    geometries}\label{4.1} 
Let us demonstrate the power of the theory developed so far by
analyzing the relation between real and complex parabolic
geometries. In particular, we will show that the theory of twistor
spaces also works in the holomorphic category.  This is of
considerable interest, since twistor correspondences usually only work
directly for some real forms of a complex semisimple Lie algebra. For
other real forms of interest, one has to restrict to real analytic
structures, pass to complexifications, and interpret the final result
on the original manifold.

Let $\frak g$ be a complex semisimple Lie algebra, $\frak
p\subset\frak g$ a parabolic subalgebra, $G$ a Lie group with Lie
algebra $\frak g$ and $P\subset G$ the parabolic subgroup
corresponding to $\frak p$. Then one can consider \textit{complex
parabolic geometries} of type $(G,P)$, defined as holomorphic
principal $P$--bundles $\Cal G\to M$ over complex manifolds $M$ which
are endowed with holomorphic Cartan connections $\om\in\Om^{1,0}(\Cal
G,\frak g)$. On the other hand, we may also consider $\frak g$ as a
real Lie algebra endowed with a parabolic subalgebra $\frak p$, and
thus consider real parabolic geometries of type $(G,P)$, which are
given by smooth principal bundles endowed with smooth Cartan
connections. See \ref{4.6} for the discussion of an interesting
example. 

If $(p:\Cal G\to M,\om)$ is a real parabolic geometry of type $(G,P)$,
then for each point $u\in\Cal G$, $\om(u)$ is a linear isomorphism
from $T_u\Cal G$ to the complex vector space $\frak g$, so $\om$
defines an almost complex structure $J^{\Cal G}$ on $\Cal G$. The
vertical subspace in $T_u\Cal G$ is just the preimage under $\om(u)$
of $\frak p$. Since this is a complex subspace, we get a complex
structure on the quotient space $T_{p(u)}M$. If we change from $u$ to
another preimage of $p(u)$, the resulting isomorphism $T_{p(u)}M\to
\frak g/\frak p$ changes by the adjoint action of an element of $P$,
which is a complex linear map. Thus, we also get a well defined almost
complex structure $J$ on $M$, and the tangent map to the projection
$p:\Cal G\to M$ is complex linear.  These almost complex structures
imply that differential forms on $M$ and $\Cal G$ with values in any
complex vector space or vector bundle split according to their complex
linearity or anti-linearity properties into $(p,q)$--types. In
particular, the adjoint tractor bundle is a complex vector bundle in
this case, so this applies to the curvature and its harmonic part
viewed as two--forms on $M$. (Note that even the real cohomologies
canonically are complex vector spaces.)

\begin{thm*}
Let $\frak g$ be a complex semisimple Lie algebra, $\frak
p\subset\frak g$ a parabolic subalgebra, $G$ a complex Lie group with
Lie algebra $\frak g$ and $P\subset G$ the parabolic subgroup
corresponding to $\frak p$. A real parabolic geometry $(p:\Cal G\to
M,\om)$ of type $(G,P)$ is actually a complex parabolic geometry
(i.e.~$M$ is a complex manifold, $p:\Cal G\to M$ a holomorphic
principal bundle and $\om$ a holomorphic Cartan connection) if and only 
if its curvature function is of type $(2,0)$. 

If the parabolic geometry $(p:\Cal G\to M,\om)$ is regular and normal,
then it suffices that the harmonic curvature is of type $(2,0)$.
\end{thm*}
\begin{proof}
We can split $L(\La^2\frak g/\frak p,\frak g)$ according to
$(p,q)$--types into $L^{2,0}\oplus L^{1,1}\oplus L^{0,2}$. Since the
adjoint action of $P$ is by complex linear maps, this is a splitting
of $P$--modules. Since $\frak g$ is a complex Lie algebra, both
$\partial$ and $\partial^*$ preserve complex multilinear maps, so for
$\Bbb E:=L^{2,0}$ we get $\square(\Bbb E)\subset\Bbb E$. One easily
verifies that $\Bbb E$ is stable under $\Bbb E$--insertions. Hence by
Corollary \ref{3.2}, in the regular normal case, the harmonic
curvature $\ka_H$ is of type $(2,0)$ if and only if the whole
curvature $\ka$ is of type $(2,0)$.

For $X\in\frak g$, we have the vector field $X^\Cal G=\om^{-1}(X)$ on
$\Cal G$, and in the proof of Proposition \ref{2.6} we have seen that
$\ka(X,Y)=[X,Y]-\om([X^\Cal G,Y^\Cal G])$. Since the bracket in $\frak
g$ is complex bilinear this formula immediately implies that the
$(0,2)$--part of $\ka$ maps $X,Y\in\frak g$ to
$-\tfrac{1}{4}\om(N(X^\Cal G,Y^\Cal G))$, where $N$ denotes the
Nijenhuis tensor of $J^\Cal G$. Consequently, the integrability of the
almost complex structure $J^\Cal G$ is equivalent to vanishing of the
$(0,2)$--component of the curvature function, and this vanishing also
implies integrability of the almost complex structure $J$ on $M$. If
the almost complex structures are integrable then the projection
$p:\Cal G\to M$ is by construction holomorphic. Equivariancy of $\om$
immediately implies that for $u\in\Cal G$ and $g\in P$ we get
$\om_{u\cdot g}\o T_ur^g=\Ad(b^{-1})\o\om_u$, which implies that
$r^g:\Cal G\to\Cal G$ is a holomorphic mapping. On the other hand,
since $\om$ reproduces the generators of fundamental vector fields,
the map $g\mapsto u\cdot g$ is holomorphic, too, and these two facts
imply that the principal right action $r:\Cal G\x P\to\Cal G$ is
holomorphic, so $p:\Cal G\to M$ is a holomorphic principal
bundle.

Given that $\Cal G$ is a complex manifold, the Cartan connection
$\om$, which by construction is a $(1,0)$--form, is holomorphic, if and
only if its exterior derivative is a $(2,0)$--form, and since the
bracket in $\frak g$ is complex bilinear, this is equivalent to the
curvature being of type $(2,0)$, which implies the result. 
\end{proof}

We can deduce from this that our results on correspondence spaces and
twistor spaces continue to hold in the realm of complex parabolic
geometries. While this is obvious for the correspondence part, it is
quite nontrivial for the twistor part. Assume that $\frak g$ is a
complex semisimple Lie algebra with standard parabolic subalgebras
$\frak q\subset\frak p\subset\frak g$ as in \ref{2.3}, $G$ is a
complex Lie group with Lie algebra $\frak g$ and $Q\subset P\subset G$
are the parabolic subgroups corresponding to $\frak q$ and $\frak
p$. If $(p:\Cal G\to N,\om)$ is a complex parabolic geometry, then
clearly $\Cal CN=\Cal G\x_PP/Q$ is a complex manifold, and $(\Cal
G\to\Cal CN,\om)$ is a complex parabolic geometry of type $(G,Q)$,
which is normal if $(\Cal G\to N,\om)$ is normal.

On the other hand, assume that $(p:\Cal G\to M,\om)$ is a complex
parabolic geometry of type $(G,Q)$, and suppose that the curvature
satisfies the restrictions of Theorem \ref{2.8} or (in the regular
normal case) \ref{3.3}. Then viewed as a real parabolic geometry, $M$
is locally isomorphic to a correspondence space $\Cal CN$ for a
(unique) real parabolic geometry of type $(G,P)$ on a manifold
$N$. But the curvature function for this parabolic geometry is the
same as the curvature function of $M$, so in particular, it has
complex bilinear values. By the theorem $N$ is a complex parabolic
geometry, and we get
\begin{kor*}
Let $(p:\Cal G\to M,\om)$ be a complex parabolic geometry of type
$(G,Q)$, which satisfies the curvature restrictions of Theorem
\ref{2.8} respectively \ref{3.3}. Then the twistor space $N$ of $M$ is
automatically a complex parabolic geometry (and in particular a
complex manifold) and the local isomorphism between $M$ and the
correspondence space $\Cal CN$ is automatically a biholomorphism.
\end{kor*}

\section{Examples and applications}\label{4}
In this section, we apply our results to three concrete
examples. These relate Lagrangian contact structures to projective
structures, certain almost CR structures to an almost complex version
of projective structures, and finally projective and almost
Grassmannian structures to a generalization of path geometries. 

\section*{Lagrangian contact structures}
In this case, the construction of the correspondence space is
described in \cite{Takeuchi} but the main result obtained there is
just that the correspondence space is locally flat if and only if the
original geometry is locally flat. Our results in this case go much
further. In one direction, we obtain a nice geometric interpretation
of the projective curvature of a torsion free affine connection, while
in the other direction we get results on contact structures and on
partial connections on projectivized cotangent bundles.

\subsection{}\label{4.2}
For $n\geq 2$ consider the Lie group $G=PSL(n+1,\Bbb R)$, the quotient
of $SL(n+1,\Bbb R)$ by its center. The Lie algebra $\frak g=\frak
s\frak l(n+1,\Bbb R)$ consists of all tracefree linear endomorphisms
of $\Bbb R^{n+1}$. Define $\frak p\subset\frak g$ to be the stabilizer
of the line through the first vector in the standard basis of $\Bbb
R^{n+1}$. Let $\frak q\subset\frak p$ be the subalgebra of those maps
which in addition preserve the hyperplane $W$ generated by the first
$n$ vectors in the standard basis. Then $\frak p$ and $\frak q$ give
rise to a $|1|$--grading, respectively a $|2|$--grading of $\frak g$,
defined by
$$
\begin{pmatrix} \frak p_0 & \frak p_1\\ \frak p_{-1} & \frak
p_0\end{pmatrix}\qquad \begin{pmatrix} \frak q_0 & \frak q^L_1 &\frak
q_2\\ \frak q_{-1}^L & \frak q_0 & \frak q^R_1\\ \frak q_{-2} & \frak
q^R_{-1} &\frak q_0\end{pmatrix},
$$
where in the first matrix, the blocks are of size $1$ and $n$, while
in the second matrix, they are of size $1$, $n-1$ and $1$. For later
use, we have indicated the finer decomposition $\frak q_{\pm 1}=\frak
q_{\pm 1}^L\oplus\frak q_{\pm 1}^R$.

The group $G$ does not act on $\Bbb R^{n+1}$ but only on the
projective space $\Bbb RP^n$. The parabolic subgroups $P,Q\subset G$
are defined as consisting of those elements whose adjoint action
preserves the filtration of $\frak g$ induced by the corresponding
subalgebra. One easily verifies that $P$ is exactly the stabilizer of
the point in $\Bbb RP^n$ corresponding to the line through the first
basis vector, so $G/P\cong\Bbb RP^n$, and $Q\subset P$ consists of
those elements which in addition stabilize the projective hyperplane
corresponding to $W$.

It is well known, see e.g. \cite{Takeuchi}, that normal parabolic
geometries of type $(G,P)$ (which are automatically regular since we
deal with a $|1|$--grading) are exactly the classical projective
structures on $n$--manifolds. Hence specifying such a geometry on an
$n$--dimensional manifold $M$ is equivalent to giving a projective
equivalence class $[\nabla]$ of torsion--free connections on the
tangent bundle $TM$. Projective equivalence means that two connections
in the class differ by the action of a one--form $\Up\in\Om^1(M)$,
i.e.~$\hat\nabla_\xi\eta=\nabla_\xi\eta+\Up(\xi)\eta+\Up(\eta)\xi$,
and it says precisely, that $\nabla$ and $\hat\nabla$ have the same
geodesics up to parametrization. Moreover, $P_0=GL(n,\Bbb R)$ and the
underlying $P_0$--bundle $\Cal G_0\to M$ is simply the full first
order frame bundle in this case.

The parabolic subalgebra $\frak q$ is of contact type, i.e.~it defines
a $|2|$--grading on $\frak g$, the components $\frak q_{\pm 2}$ have
dimension one, and the bracket $\frak q_{-1}\x\frak q_{-1}\to\frak
q_{-2}$ is non degenerate. Therefore we get an example of a
\textit{parabolic contact structure}, i.e.~regular normal parabolic
geometries of type $(G,Q)$ have an underlying contact structure.
These geometries have been first studied in \cite{Takeuchi}, where
they were called Lagrangian contact structures (although Legendrian
contact structures would be closer to the usual terminology). To
describe them, consider a contact structure $H\subset TM$ on a smooth
manifold $M$ of dimension $2n+1$. This means that $H$ is a subbundle
of rank $2n$ which is maximally non--integrable in the sense that the
skew symmetric bundle map $\Cal L:H\x H\to TM/H$ induced by the Lie
bracket of vector fields is non--degenerate. A subbundle $E\subset H$
is called \textit{isotropic} if the restriction of $\Cal L$ to $E\x E$
vanishes identically. Elementary linear algebra shows that isotropic
subbundles have rank at most $n$, and isotropic subbundles of maximal
rank are usually called \textit{Legendrian} or \textit{Lagrangian}.

\subsection*{Definition} A \textit{Lagrangian contact structure} on a
smooth manifold $M$ of dimension $2n+1$ is a contact structure
$H\subset TM$ together with a fixed splitting $H=H^L\oplus H^R$ of the
contact subbundle as a direct sum of two rank $n$ isotropic
subbundles. 

\bigskip

Note that if $(M,H=H^L\oplus H^R)$ is a Lagrangian contact structure,
then $\Cal L$ induces an isomorphism $H^R\cong L(H^L,TM/H)$, so the
two subbundles are almost dual to each other. 

In \cite{Takeuchi}, the author used Tanaka's prolongation procedure
from \cite{Tanaka79} to construct a parabolic geometry from a
Lagrangian contact structure. We sketch a construction using the
procedure from \cite{Cap-Schichl}, which uses a simpler description of
the underlying structures. We have noted above that the Lie bracket
$[\ ,\ ]:\frak q_{-1}\x\frak q_{-1}\to\frak q_{-2}$ is
non--degenerate, and one observes that the subspaces $\frak q_{-1}^L$
and $\frak q_{-1}^R$ are isotropic. Next, one verifies that the
adjoint action identifies the subgroup $Q_0\subset Q$ with the group
of all linear isomorphisms $\ph:\frak q_{-1}\to\frak q_{-1}$ which
preserve the decomposition $\frak q_{-1}^L\oplus\frak q_{-1}^R$ and
have the property that there exists a linear isomorphism $\ps:\frak
q_{-2}\to\frak q_{-2}$ such that $[\ph(X),\ph(Y)]=\ps([X,Y])$ for all
$X,Y\in\frak g_{-1}$. Of course, the map $\ps$ is then uniquely
determined by $\ph$. In the description of \cite{Cap-Schichl} for a
parabolic geometry of type $(G,Q)$ one first needs a manifold of
dimension $2n+1$ together with a subbundle $H\subset TM$ of rank
$\dim(\frak g_{-1})=2n$. Then regular normal parabolic geometries of
type $(G,Q)$ are in bijective correspondence with reductions of the
associated graded $H\oplus TM/H$ to the structure group $Q_0$ such
that the bundle map $\La^2H\to TM/H$ obtained from the Lie bracket of
vector fields looks like the bracket $[\ ,\ ]:\La^2\frak
q_{-1}\to\frak q_{-2}$ in each fiber. The condition on the bracket
exactly means that $H$ defines a contact structure on $M$, and the
reduction to $Q_0$ is exactly equivalent to the decomposition
$H=H^L\oplus H^R$ into a sum of Lagrangian subbundles.

Now we can apply the construction of correspondence spaces to this
case, which recovers all the results of \cite{Takeuchi}:
\begin{thm*}
Let $[\nabla]$ be a projective structure on a smooth manifold $N$, and
let $(p:\Cal G\to N,\om)$ be the associated parabolic geometry of type
$(G,P)$. Then the correspondence space $\Cal CN$ for $\frak
q\subset\frak p$ is the projectivized cotangent bundle $\Cal
P(T^*N)$. The induced Lagrangian contact structure $H=H^L\oplus H^R$
on $\Cal P(T^*N)$ has the following form: $H\subset T\Cal P(T^*N)$ is
the canonical contact structure, $H^R$ is the vertical subbundle of
$\Cal CN\to N$, and $H^L$ is obtained from the horizontal
distributions of the connections in the projective class.
\end{thm*}
\begin{proof}
Since the Killing form induces a duality between $\frak g/\frak p$ and
$\frak p_1$, the cotangent bundle $T^*N$ is the associated bundle
$\Cal G\x_P\frak p_1$. By definition, the subgroup $Q\subset P$
preserves the line $\frak q_2\subset\frak p_1$, and one immediately
verifies that this property characterizes $Q$. Passing to the
projectivization, we see that $\Cal P(\frak p_1)\cong P/Q$, which
implies that $\Cal CN=\Cal G/Q=\Cal G\x_P(P/Q)$ is exactly the
projectivized cotangent bundle $\Cal P(T^*N)$. 

The tangent bundle to $\Cal P(T^*N)$ is $\Cal G\x_Q(\frak g/\frak q)$,
while the tangent bundle to $N$ is $\Cal G\x_P(\frak g/\frak p)$. The
tangent map of the projection $\pi:\Cal P(T^*N)\to N$ corresponds in
this picture exactly to the projection $\frak g/\frak q\to \frak
g/\frak p$. The contact distribution $H\subset T\Cal P(T^*N)$ is given
by $\Cal G\x_Q((\frak q_{-1}\oplus\dots\oplus\frak q_2)/\frak
q)$. Since $\frak q_{-1}\oplus\dots\oplus\frak q_2$ is exactly the
annihilator of $\frak q_2$ with respect to the Killing form, we
conclude that for a point $\ell\in\Cal P(T^*N)$ (i.e.~$\ell$ is a line
in $T^*_{\pi(\ell)}M$), we have $H_\ell=\{\xi:\ell(T\pi\cdot\xi)=0\}$,
so we exactly get the canonical contact structure on $\Cal
P(T^*N)$. Since $\frak q_{-1}^R=\frak q\cap\frak p$, the subbundle
$H^R$ consists of those tangent vectors, which project to zero, so
this is exactly the vertical subbundle in $T\Cal P(T^*N)$.

Any connection $\nabla$ in the projective class induces a linear
connection on $T^*M$, which gives rise to a vertical projection from
$TT^*M$ onto the vertical subbundle. This vertical projection is
characterized by the fact that the covariant derivative of a one form
$\ph$ is the nontrivial component of the composition of the vertical
projection with the tangent map $T\ph$. Factoring to the
projectivization we see, that any linear connection on $TN$ gives rise
to a vertical projection on $T\Cal P(T^*N)$. For projectively
equivalent connections $\nabla$ and $\hat\nabla$, one easily computes
that on the level of one--forms one gets
$\hat\nabla_\xi\ph=\nabla_\xi\ph-\Up(\xi)\ph-\ph(\xi)\Up$. Passing
from $T_\ph T^*N$ to $T_{[\ph]}\Cal P(T^*N)$ means exactly factoring
by the line generated by $\ph$ in the vertical subspace, and the
subbundle $H$ is characterized by $\xi\in H_{[\ph]}$ iff $\ph(\xi)=0$,
which shows that the restriction of the vertical projection to
$H\subset T\Cal P(T^*N)$ depends only on the projective class. One
verifies directly that this construction really describes the bundle
$H^L$.
\end{proof}

Note that the subbundle $H^L\subset H$ which is complementary to the
vertical subbundle, is very similar to a connection on the fiber
bundle $\pi:\Cal P(T^*N)\to N$. The difference to a true connection is
that the vertical projection is only defined on the subbundle
$H\subset T\Cal P(T^*N)$. Thus, $H^L$ defines a \textit{partial
connection} on $\Cal P(T^*N)$.

\subsection{Harmonic curvature}\label{4.3}
To proceed further, we have to describe the harmonic curvature both
for projective structures and for Lagrangian contact structures. The
structure of this curvature is essentially different for $n=2$ and
$n\geq 3$. We discuss the case $n\geq 3$ in detail and make some
remarks on the case $n=2$ below.

We have to compute the second cohomologies $H^2(\frak p_-,\frak g)$
and $H^2(\frak q_-,\frak g)$, and this is a simple exercise in
applying Kostant's version of the Bott--Borel--Weil theorem, see
\cite{Kostant}. Since we are dealing with the split real form $\frak
s\frak l(n+1,\Bbb R)$ of $\frak s\frak l(n+1,\Bbb C)$ here, the real
cohomologies look exactly as the complex cohomologies. Using the
Dynkin--diagram notation from chapter 3 of \cite{Baston-Eastwood},
i.e.~the numbers over the nodes indicate the coefficient of the
corresponding fundamental weight in the highest weight of the dual
representation, and the algorithms from section 8.5 of
\cite{Baston-Eastwood} one gets $H^2(\frak p_-,\frak
g)=\xbb{-4}{1}{2}$ respectively $\xbbdb{-4}{1}{1}{0}{1}$ for the
projective case. This is immediately seen to be exactly the
irreducible component of highest weight in $\La^2(\frak
p_{-1})^*\otimes L(\frak p_{-1},\frak p_{-1})$. Following the
description in \ref{2.2} one shows that the harmonic curvature is
represented by trace free part of the curvature of any connection in
the projective class.

For the Lagrangian contact structures, the situation is a little more
complicated, since there are three irreducible components in the
curvature. One obtains the following picture:

\smallskip

\begin{center}
\begin{tabular}{|l|l|}
\hline
\parbox[c][1.2\totalheight][c]{68pt}{$\xbbdx{-4}{1}{1}{0}{1}$} & torsion $T^R:\La^2H^R\to H^L$\\
\hline
\parbox[c][1.2\totalheight][c]{68pt}{$\xdbbx{1}{0}{1}{1}{-4}$} & torsion $T^L:\La^2H^L\to H^R$\\
\hline
\parbox[c][1.2\totalheight][c]{68pt}{$\xbdbx{-3}{2}{0}{2}{-3}$} &
curvature $\rho:H^L\x H^R\to L(H^L,H^L)$\\ \hline
\end{tabular}
\end{center}

\smallskip

The first column contains the representations for $n\geq 4$, and the
second column contains the corresponding geometric object. For $n=3$,
the representations are slightly different, but the geometric objects
remain the same. The representations are computed as before and the
geometric objects are determined by identifying them with components
in $\La^2(\frak g/\frak q,\frak g)$ and passing to the corresponding
bundles. The torsions $T^L$ and $T^R$ are easy to interpret
geometrically. Since $H^L$ is isotropic, the bracket of two sections
of $H^L$ is a section of $H$ and thus can be projected to
$H^R$. Similarly as in Proposition 3.4 of \cite{Cap-Schmalz} one shows
that the resulting tensor field coincides with $T^L$ up to a nonzero
multiple. In particular, $T^L$ is exactly the obstruction to
integrability of the subbundle $H^L\subset TM$. The interpretation of
$T^R$ is completely analogous. To compute the remaining harmonic
curvature component $\rho$, one has to choose a Weyl--structure and
compute the appropriate harmonic part of its curvature.

Since $\frak q_{-1}\cap\frak p=\frak q_{-1}^R$, the distribution
corresponding to $\frak p/\frak q\subset\frak g/\frak q$ is exactly
$H^R$. Integrability of this subbundle is equivalent to vanishing of
$T^R$. Moreover, Theorem \ref{3.3} implies that local leaf spaces
carry an induced projective structure (an hence the original structure
is locally isomorphic to a correspondence space) if and only also
$\rho$ vanishes identically. In particular, for a correspondence space
the torsion $T^L$ is a complete obstruction against local flatness,
and we obtain:
\begin{thm*}
The projective Weyl curvature of a projective structure $[\nabla]$ on
a smooth manifold $N$ of dimension $n\geq 3$ is exactly the
obstruction to integrability of the bundle $H^L$ on $\Cal P(T^*N)$,
i.e.~to flatness of the induced partial connection on $\Cal P(T^*N)\to
N$.
\end{thm*}
Having a twistor correspondence in the classical sense would mean that
one starts with a projective manifold $(N,[\nabla])$, forms the
correspondence space $\Cal CN$ and then a twistor space with respect
to the parabolic $\tilde{\frak p}\supset\frak q$, which is the
stabilizer of the line spanned by the last vector in the standard
basis of $\Bbb R^{n+1}$. The corresponding distribution is $H^L$, so
the theorem shows that this is possible only in the locally
projectively flat case, in which one recovers projective duality.

\subsection{}\label{4.4}
Let us next look at a Lagrangian contact structure $(M,H^L\oplus H^R)$
which admits a twistor space but not a projective structure on this
twistor space. So assume that $H^R$ is integrable and let $\ps:U\to N$
be a local leaf space, see \ref{2.6}. For any $x\in U$, we have the
subspace $H_x\subset T_xU$, and by construction the image of this
subspace is a hyperplane in $T_{\ps(x)}N$. The annihilator of this
hyperplane is a line in $T^*_xN$, and thus a point
$\tilde\ps(x)\in\Cal P(T^*N)$. Clearly, this defines a smooth mapping
$\tilde\ps:U\to\Cal P(T^*N)$ such that $\pi\o\tilde\ps=\ps$, where
$\pi:\Cal P(T^*N)\to N$ is the projection. To understand $\tilde\ps$,
it is better to view $\Cal P(T^*N)$ as the space $Gr_n(TN)$ of
hyperplanes in the tangent space of $N$. By construction, the kernel
of $T_x\tilde\ps$ has to be contained in $\ker(T_x\ps)=H^R_x$. For
$\xi\in H^R_x$ the image $T_x\tilde\ps\cdot\xi$ is contained in the
vertical subspace $V_{\tilde\ps(x)}Gr_n(TN)$. This is the tangent
space to the fiber at $\tilde\ps(x)$, so it can be identified with the
space $L(\tilde\ps(x),T_{\ps(x)}N/\tilde\ps(x))$ of linear maps. Via
$T_x\ps$, this space is isomorphic to $L(H_x/H^R_x,T_xM/H_x)$. Going
through the identifications one sees that
$T_x\tilde\ps\cdot\xi(x)(\eta)=\Cal L(\xi(x),\eta)$ for all $\eta\in
H_x$, so non degeneracy of $\Cal L$ implies that $T_x\tilde\ps$ is
injective.

Since $M$ and $\Cal P(T^*N)$ have the same dimension, we conclude
that, possibly shrinking $U$, that $\tilde\ps$ is a diffeomorphism
from $U$ onto an open subset $V$ of $\Cal P(T^*N)$. Moreover, from the
construction it is obvious that this is a contact diffeomorphism, and
the integrable Lagrangian subbundle $H^R$ is exactly mapped to the
vertical subbundle of $\Cal P(T^*N)$. Notice that the complementary
Lagrangian subbundle $H^L$ was not used at all in the
construction. Hence we obtain the following strengthening of the
Darboux theorem:
\begin{thm*}
Let $(M,H)$ be a contact manifold of dimension $2n+1$, and let
$H^R\subset H$ be an integrable Lagrangian subbundle. Then locally $M$
is contact diffeomorphic to the projectivized tangent bundle of $\Bbb
R^{n+1}$ in such a way that $H^R$ is mapped to the vertical subbundle.
\end{thm*}

\subsection{}\label{4.5}
We now know that a Lagrangian contact manifold $(M,H^L\oplus H^R)$
which has the property that $T^R$ is identically zero is locally
contact diffeomorphic to an open subset of the projectivized tangent
bundle of its twistor space in such a way that $H^R$ is mapped to the
vertical subbundle. On such an open subset the image of $H^L$ defines
a complement to the vertical subbundle in the contact subbundle, so we
can view this as a locally defined partial connection.

As we have observed, any linear connection $\nabla$ on $TN$ gives rise
to such a partial connection (even globally defined), and moreover,
this partial connection depends only on the projective class of
$\nabla$. By Theorem \ref{3.3}, a locally defined partial connection
comes from a linear connection on $TN$ if and only if the curvature
$\rho\in\Ga((H^L)^*\otimes(H^R)^*\otimes L(H^L,H^L))$ of the
corresponding Lagrangian contact structure is identically zero, so we
get:
\begin{thm*}
Let $N$ be a smooth manifold of dimension $n\geq 3$ and let
$U\subset\Cal P(T^*N)$ be an open subset. Let $H$ be the restriction
of the contact subbundle to $U$ and let $H^L\subset H$ be a Lagrangian 
subbundle, which is complementary to the vertical subbundle $H^R$. Then
$H^L$ is obtained from a linear connection on $TN$ as described in
the end of \ref{4.2} if and only if the curvature $\rho$ of the
Lagrangian contact structure $(U,H^L\oplus H^R)$ vanishes identically. 
\end{thm*}

\subsection{The case $n=2$}\label{4.5a} 
In this case, the structure of the harmonic curvature is completely
different. Rather than having two torsions and one curvature, one has
two curvatures in this case. The Lagrangian contact structures in this
case admit a natural interpretation as path geometries, so this can
also be viewed as a special case of the structures discussed in
\ref{4.7}. Also projective structures in two dimensions behave
differently than in general. For dimensional reasons the tracefree
part of the curvature of any linear connection in two dimensions
vanishes identically. The complete obstruction to local projective
flatness is a tensor analogous to the Cotton--York tensor in conformal
geometry. Passing to the correspondence space, this corresponds
exactly to one of the two curvatures. Since the bundles $H^L$ and
$H^R$ have rank one, they are always integrable. In particular, in
this special dimension, one can obtain a twistor correspondence in the
classical sens in non--flat situations. This looks particularly nice
in the picture of path geometries, viewed as encoding second order
ODE's, compare with \ref{4.7}. In this picture, the twistor space is
the space of all solutions of the equation. The holomorphic version of
this correspondence was used in \cite{Hitchin} to study Schlesinger's
equation.

\section*{Elliptic partially integrable almost CR manifolds}

\subsection{}\label{4.6}
This is an almost complex version of the situation discussed in
\ref{4.5a} above. We consider $G=PSL(3,\Bbb C)$ as a real Lie group,
the Borel subgroup $Q\subset G$ (i.e.~$Q$ is the stabilizer of the
standard flag $\Bbb C\subset\Bbb C^2\subset\Bbb C^3$) and the
stabilizer $P$ of the complex line $\Bbb C\subset\Bbb C^3$. This case
is a bit involved and will be taken up in detail elsewhere. Here we
only give a brief outline. From their definition, real parabolic
geometries of type $(G,Q)$ are six dimensional smooth manifolds
endowed with an almost complex structure and two complementary complex
line bundles $H^+,H^-\subset TM$, such that, with $H=H^+\oplus H^-$,
the tensorial map $\Cal L:H\otimes H\to TM/H$ induced by the Lie
bracket is complex bilinear and non--degenerate. Building on earlier
work in \cite{Schmalz-Slovak} it has been shown in \cite{Cap-Schmalz}
that flipping the almost complex structure on the subbundle $H^+$
leads to an equivalence of categories between the category of regular
normal parabolic geometries of type $(G,Q)$ and the category of
elliptic partially integrable almost CR manifolds of CR dimension and
codimension two. Since we are dealing with the underlying real Lie
algebra of a complex Lie algebra, the structure of torsions and
curvatures is rather complicated. In the notation of
\cite{Cap-Schmalz} there are the following irreducible components of
the harmonic curvature:

\smallskip

\begin{center}
\begin{tabular}{|l|l|l|l|}
\hline
$S^\pm$ & $(0,2)$ & $TM/HM\x H^{\pm} M\to TM/HM$ & almost complex structure\\
\hline
$N^\pm$ & $(1,1)$ & $H^+M\x H^-M\to H^\pm M$ & almost CR structure\\
\hline
$T^\pm$ & $(1,1)$ & $\La^2H^\pm M\to H^\mp M$ & integrability of $H^\pm M$\\
\hline
$\rho^\pm$ & $(2,0)$ & $TM/HM\x H^{\pm}\to (H^{\pm})^*$ & \\
\hline
\end{tabular}
\end{center}

\smallskip

The second column contains the $(p,q)$--types of the components, while
the last column indicates the structure whose integrability is
obstructed by the component. The first three lines correspond to the
components of torsion type, while the two components in the last line
are curvatures. From this table and Theorem \ref{4.1} we immediately
conclude that complex regular normal parabolic geometries of type
$(G,Q)$ are characterized among the real ones by vanishing of $S^\pm$,
$N^\pm$, and $T^\pm$, which is easily seen to be equivalent to torsion
freeness. For such a complex parabolic geometry the almost CR
structure is automatically integrable (since $N^\pm$ vanishes) and
real analytic, which implies that torsion free elliptic CR manifolds
are always locally embeddable, see \cite{Cap-Schmalz}. 

On the other hand, it turns out that normal parabolic geometries of
type $(G,P)$ are an almost complex analog of two dimensional
projective structures. Given an almost complex manifold $(\Cal N,J)$
we call two linear connections $\nabla$ and $\hat\nabla$
\textit{projectively equivalent} if there exists a smooth
$(1,0)$--form $\Up$ on $\Cal N$ such that
$\hat\nabla_\xi\eta=\nabla_\xi\eta+\Up(\xi)\eta+\Up(\eta)\xi$. Note
that does \textit{not} imply projective equivalence in the real sense,
since complex multiples of $\xi$ and $\eta$ are involved. One easily
shows projectively equivalent connections have the same torsion and if
$\nabla J=0$ for some connection $\nabla$, then the same is true for
any projectively equivalent connection. Now we define a
\textit{compatible projective structure} on $(\Cal N,J)$ to be a
projective class $[\nabla]$ of connections such that $\nabla J=0$ and
the torsion is of type $(0,2)$. This is the best possible
normalization of the torsion and it implies that the torsion is given
by $-1/4$ times the Nijenhuis tensor of $J$. Similar to the case of
classical projective structures one proves that normal parabolic
geometries of type $(G,P)$ are exactly compatible projective
structures on almost complex manifolds $(\Cal N,J)$ of real dimension
four.

The harmonic curvature of these compatible projective structures
consists of three irreducible components, one in each of the types
$(0,2)$, $(1,1)$, and $(2,0)$. The $(0,2)$--part is a multiple of the
Nijenhuis tensor of $J$, the $(1,1)$--part is essentially the
obstruction against the projective class locally containing
holomorphic connections, while the $(2,0)$--part is exactly the
complex analog of the projective curvature in two dimensions.

Compatible projective structures are by far simpler in nature than
partially integrable almost CR structures. In particular, they exist
on any almost complex manifold $(\Cal N,J)$ and the set of all such
structures is very easy to describe. Hence in this case already the
correspondence space construction leads to a very interesting result:
\begin{thm*}
Let $(\Cal N,J)$ be an almost complex manifold of real dimension
four. Then any choice of a compatible projective structure $[\nabla]$
on $\Cal N$ endows the space $M:=\Cal P_{\Bbb C}(T\Cal N)$ of complex
lines in the tangent spaces of $\Cal N$ with an elliptic CR structure
of CR dimension and codimension 2.

For this CR structure, the components $T^-$, $N^\pm$, $S^-$ and
$\rho^-$ of the harmonic curvature vanish identically, while the three
remaining components $S^+$, $T^+$, and $\rho^+$ correspond directly to
the three components of the harmonic curvature of $(\Cal
N,J,[\nabla])$ of the respective $(p,q)$--type. The CR automorphism
group of this structure coincides with the group of projective
automorphisms of $(\Cal N,J,[\nabla])$. 
\end{thm*}

On the other hand, for a general regular normal parabolic geometry
of type $(G,P)$, one can use the tools developed in this paper to
give a twistorial interpretation of the components of the harmonic
curvature. The subbundle in $TM$ corresponding to $\frak p\supset\frak
q$ is the bundle $H^-$. Here there are three steps:

\smallskip

\begin{center}
\begin{tabular}{|l|l|}
\hline
existence of a local leaf space $\Cal N$ & $T^-=0$ \\
\hline
almost complex structure $J$ on $\Cal N$ & $N^+=0,S^-=0$ \\
\hline
compatible projective structure on $(\Cal N,J)$ & $N^-=0,\rho^-=0$\\
\hline
\end{tabular}
\end{center}

\smallskip 

In the second step, one asks whether the almost complex structure on
$M$ descends to a local leaf space $\Cal N$, and verifying that this
is equivalent to the vanishing of $N^+$ and $S^-$ needs a bit of extra
work. If this is the case, then $M$ is locally diffeomorphic to the
complex projectivization of $T\Cal N$. Under this condition, the
geometry on $M$ can then be described locally as an almost complex
version of a path geometry, or equivalently as a partial connection on
$\Cal P_{\Bbb C}(T\Cal N)$. The last step directly follows from
Theorem \ref{3.3}.

Let us finally mention that in this case there is the possibility to
form a twistor correspondence in the classical sense. Given a
compatible projective structure $(\Cal N,J,[\nabla])$ one can form the
correspondence space $M$, and the harmonic curvature is encoded in
$S^+$, $T^+$, and $\rho^+$. If the component corresponding to $T^+$
vanishes, then the bundle $H^+$ is integrable, and one can form a
local leaf space $Z$. If one wants this leaf space to carry an induced
almost complex structure, then also the component corresponding to
$S^+$ has to vanish, and one is in the holomorphic situation as
discussed in \cite{Hitchin}. 

\section*{Almost Grassmannian structures}

This example, which we only outline briefly, contains the twistor
theory for paraconformal manifolds of \cite{Bailey-Eastwood} as well
as twistor theory for conformal four manifolds in split signature
case. On the other hand, via the notion of path geometries, there is a
relation to the geometric theory of systems of second order ODE's,
see \cite{Grossman}.

\subsection{}\label{4.7}
For $n\geq 2$ we consider $G=PSL(n+2,\Bbb R)$ and the parabolic
subgroups $P$ and $\hat P$ of $G$ defined as the stabilizer of a point
respectively a line containing that point in $\Bbb RP^{n+1}$. Then
also $Q:=P\cap \hat P$ is a parabolic subgroup of $G$. From \ref{4.2}
we know that normal parabolic geometries of type $(G,P)$ are classical
projective structures on $n+1$--dimensional manifolds. On the other
hand, normal parabolic geometries of type $(G,\hat P)$ are exactly
almost Grassmannian (also called paraconformal) structures, see
\cite{Bailey-Eastwood}. Essentially they are defined as smooth
manifolds $\hat N$ of dimension $2n$ together with an isomorphism
$T\hat N\cong E^*\otimes F$, where $E$ and $F$ are auxiliary bundles
on $\hat N$ of rank $2$ and $n$, respectively. In the special case
$n=2$, one obtains exactly four dimensional split signature conformal
manifolds. It should also be noted that for a different real form, one
obtains almost quaternionic structures.

We have described the $|1|$--grading corresponding to $\frak
p\subset\frak g=\frak s\frak l(n+2,\Bbb R)$ in \ref{4.2}. For
$\hat{\frak p}$, one also obtains a $|1|$--grading, which has the same
form, except that the blocks have size $2$ and $n$ rather than $1$ and
$n+1$. For $\frak q=\frak p\cap\hat{\frak p}$, one obtains a
$|2|$--grading on $\frak g$, with $\frak q_{\pm 2}=\frak
p_\pm\cap\hat{\frak p}_\pm$ of dimension $n$, and $\frak q_{\pm
1}=\frak q_{\pm 1}^L\oplus \frak q_{\pm 1}^R$. Here $\frak q_{\pm
1}^L=\frak q_{\pm 1}\cap\hat{\frak p}_\pm$ is one--dimensional while
$\frak q_{\pm 1}^R=\frak q_{\pm 1}\cap\frak p_\pm$ has dimension $n$.

Using the main result of \cite{Cap-Schichl}, one shows that regular
normal parabolic geometries of type $(G,Q)$ are $2n+1$--dimensional
manifolds $M$ equipped with complementary subbundles $H^L,H^R\subset
TM$ of rank $1$ and $n$, respectively, which have the property that
the bracket of two sections of $H^R$ is a section of $H:=H^L\oplus
H^R$ and the tensorial map $\Cal L:H\x H\to TM/H$ induced by the Lie
bracket is non--degenerate. Note that these assumptions in particular
imply that $H^R\cong L(H^L,TM/H)$.

From \ref{4.2} we know that the harmonic curvature for projective
structures is the Weyl curvature, i.e.~the tracefree part of the
curvature of any connection in the projective class. The harmonic
curvatures for almost Grassmannian structures are also well
known. There are always two irreducible component, but there is an
important difference between the cases $n=2$ and $n\geq 3$. For $n=2$
both of the two components are curvatures, and they correspond to the
self dual and the anti self dual part of the Weyl curvature for
four dimensional conformal structures. On the other hand, for $n\geq
3$, one of the two components is a torsion, while the other one is a
curvature. The torsion is exactly the obstruction to the existence of
a torsion free connection preserving the almost Grassmannian
structure, i.e.~to the structure being Grassmannian.

It turns out that the harmonic curvature for geometries of type
$(G,Q)$ always has three irreducible components and the highest
weights of the corresponding representations are exactly the
restrictions to $\frak q_0=\frak p_0\cap\hat{\frak p}_0$ of the
highest weights of the representations corresponding to the three
curvature components discussed above. The component corresponding to
the projective Weyl curvature is realized by a torsion $T:H^L\x
TM/H\to H^R$, while the component corresponding to the curvature on
the almost Grassmannian side is realized by a curvature $\rho:H^R\x
TM/H\to L(H^R,H^R)$. For $n=2$ the second curvature on the
Grassmannian side is realized by a torsion $\tau:\La^2H^R\to H^L$,
while for $n\geq 3$ the torsion $\tau$ corresponding to the torsion on
the Grassmannian side has the form $\tau:\La^2H^R\to TM/H$. Hence for
$n\geq 3$ the torsion $\tau$ is homogeneous of degree zero, so it has
to vanish for regular normal parabolic geometries. In particular,
correspondence spaces of almost Grassmannian structures with nontrivial
torsion are examples of parabolic geometries of type $(G,Q)$ which are
normal but not regular.

Similarly as in \ref{4.2} one verifies that for a projective manifold
$(N,[\nabla])$ of dimension $n+1$ the correspondence space for $\frak
q\subset\frak p$ is the projectivized tangent bundle $\Cal P(TN)$. The
bundle $H^R$ is the vertical subbundle of $\Cal P(TN)\to N$ and $H$ is
the tautological bundle, whose fiber at a point $\ell$ consists of all
tangent vectors whose projection to $TN$ lie in the line $\ell$. The
line bundle $H^L\subset H$ which is complementary to $H^R$ is
constructed from the horizontal lifts of the connections in the
projective class as in \ref{4.2}. For an almost Grassmannian manifold
$\hat N$, the correspondence space for $\frak q\subset\hat{\frak p}$
is the projectivization of the auxiliary rank two bundle $E$, $H^L$
is the vertical bundle, and $H^R$ is constructed from the almost
Grassmannian structure.

For $n=2$, Theorem \ref{3.3} implies that integrability of $H^R$ is
equivalent to vanishing of $\tau$. For $n\geq 3$ one shows using Lemma
3.2 of \cite{Cap-Schmalz} that $H^R$ is integrable for any regular
normal parabolic geometry of type $(G,Q)$. If this is satisfied and
$N$ is a local leaf space, then similarly as in \ref{4.4} one shows
that $M$ is locally diffeomorphic to $\Cal P(TN)$ and this
diffeomorphism maps the subbundles $H^R\subset H\subset TM$ to the
vertical respectively the tautological subbundle of $T\Cal P(TN)$ as
described above. The subbundle $H^L$ is then mapped to a line
subbundle in $T\Cal P(TN)$ which is complementary to the vertical
subbundle in the tautological subbundle and contains the complete
information about the local geometry on $M$. Conversely, the
non--integrability properties of the tautological subbundle in $T\Cal
P(TN)$ imply that such a complementary line bundle always gives rise
to a regular normal parabolic geometry of type $(G,Q)$. These
complementary line bundles are exactly the \textit{path--geometries}
as defined for example in \cite{Grossman} (via differential systems),
and using Theorem \ref{3.3} we obtain:
\begin{thm*}
(1) Any path geometry on a smooth manifold $N$ of dimension $n+1$
gives rise to a regular normal parabolic geometry of type $(G,Q)$. If
$n=2$ then the torsion $\tau$ of this geometry vanishes
identically. Conversely, if $n\geq 3$ or $\tau=0$ any such parabolic
geometry on a manifold $M$ locally admits a twistor space $N$
corresponding to $P\supset Q$, and is locally isomorphic to a path
geometry on this twistor space.

\noindent
(2) The path geometry in (1) comes from a projective structure on $N$
if and only if $\rho$ vanishes identically, and then $T$ corresponds to
the projective Weyl curvature.

\noindent
(3) For any regular normal parabolic geometry of type $(G,Q)$ on $M$,
there exists a local twistor space $\hat N$ corresponding to $\hat
P\supset Q$. The structure on $M$ descends to a Grassmannian
(respectively anti self dual conformal) structure on $\hat N$ if and
only if $T$ vanishes identically.
\end{thm*}

Via the correspondence between path geometries and systems of second
order ODE's, part (2) describes when such a system can be written
as the geodesic equation for some connection. The definition of
torsion free path geometries in \cite{Grossman} is easily seen to be
equivalent to vanishing of the torsion $T$. Parts (1) and (2) imply
that a projective structure leading to a torsion free path geometry is
locally flat, which is Theorem 1 of \cite{Grossman}. Part (3) implies
that for torsion free equations the structure descends to a
Grassmannian (respectively anti self dual conformal) structure on
$\hat N$ (which is constructed in \cite{Grossman} via a Segre
structure). Hence the curvature descends to $\hat N$, which
generically leads to explicit solutions for equations corresponding to
torsion free path geometries.

Applying part (1) to the correspondence space of an almost
Grassmannian (respectively split signature conformal) shows that the
conditions for existence of a twistor space are vanishing torsion
respectively anti self duality.  This recovers the standard twistor
theory for these structures, see \cite{Bailey-Eastwood} for the
Grassmannian case. From (2) we conclude that a local geometric
structure on the twistor space is only available in the locally flat
case.

\end{document}